\documentclass{article}

\usepackage{Jared-MK}
\usepackage{graphicx} 
\usepackage{mathtools} % http://ctan.org/pkg/mathtools

\newcommand{\bg}{\overline{g}}
\newcommand{\bh}{\overline{h}}

\newcommand{\eps}{\epsilon}

\newcommand{\Leu}{\eps^2 \Delta_g - W''(u)}
\newcommand{\Rwe}{R_{\omega, \eps}}
\newcommand{\bw}{\overline{w}}

\newtheorem*{theorem*}{Theorem}

\usepackage{thmtools}
\usepackage{thm-restate}
\usepackage{cleveref}

\title{A Dirichlet-to-Neumann Map for the Allen-Cahn Equation on Manifolds with Boundary}
\author{Jared Marx-Kuo}

\date{April 2023}
\begin{document}

\maketitle

\begin{abstract}
\noindent We study the asymptotic behavior of Dirichlet minimizers to the Allen--Cahn equation on manifolds with boundary, and we relate the Neumann data to the geometry of the boundary. We show that Dirichlet minimizers are asymptotically local in orders of $\eps$ and compute expansions of the solution to high order. A key tool is showing that the linearized allen-cahn operator about is invertible at the heteroclinic solution, on functions with $0$ boundary condition. We apply our results to separating hypersurfaces in closed Riemannian manifolds. This gives a projection theorem about Allen--Cahn solutions near minimal surfaces, as constructed by Pacard--Ritore.
%\noindent We consider $Y^{n-1} \subseteq (M^{n}, g)$, a two sided, separating hypersurface, along with energy minimizers, $u_{\eps, Y}$ of the Allen--Cahn energy with Dirichlet conditions on $Y$. We show that for $Y$ a $C^{3,\alpha}$ surface with bounded geometry, $\partial_{\nu^{\pm}} u^{\pm}_{\eps, Y} = (\eps \sqrt{2})^{-1} \pm \sigma_0 H_Y + O(\eps^{1-\alpha})$. In general, if $Y$ is $C^{k+3, \alpha}$, $\partial_{\nu}u_{\eps, Y}$ can be expanded to order $\eps^k$ with computable coefficients. We use this to prove that the projection of Allen--Cahn solutions constructed near minimal surfaces in Pacard--Ritore reflects the perturbation between the minimal surface and the zero set of the solution. Along the way, we demonstrate a better Schauder estimate for the linearized Allen--Cahn operator, $L = \eps^2\Delta_{g} - W''(g_{\eps})$, acting on functions with Dirichlet condition $\phi\Big|_Y \equiv 0$. While we are motivated by two sided separating hypersurfaces, our techniques work for manifolds with boundary $(M^n, g)$, such that $\partial M = Y$. 
\end{abstract}

\tableofcontents 

\begin{comment}
\section{To Do}
\begin{itemize}

\item (\textcolor{blue}{done}) Make $g(t)$ the true heteroclinic and $\bg$ the modified one 
\item (\textcolor{blue}{done}) Need to reframe 
\[
||\phi||_{C^{2,\alpha}_{\eps}(M^+)} \leq K ||\Leu(\phi)||_{C^{\alpha}_{\eps}(M^+)}
\]
And so need to change all of the bounds from $Y \times [0, - \omega \eps \ln(\eps)) \to M^+$

\item (\textcolor{blue}{done}) If I can prove this, then we should have that $\Leu$ has no kernel in $C^{2,\alpha}_{\eps, 0}$ - need to make this into a proposition!

\item (\textcolor{blue}{done})	In particular, need to let 
\[
g_{\eps}(t) = [1 - \chi(t/(-\omega \eps\ln(\eps)))] g(t) + \chi(t/ (-\omega \eps\ln(\eps)))
\]
where $\chi(t)$ is $0$ on $t < 1$, $1$ on $t > 2$, and goes from $0 \to 1$ on $[1,2]$. This tells us that 
\[
[\eps^2 \partial_t^2 - W''(g_{\eps})] g_{\eps} = O(\eps^{\omega})
\]

\item (\textcolor{blue}{done}) Yeah check for all instances of ``$\ln(\eps)$"

\item (\textcolor{blue}{done}) Coefficient $\phi_t(s,0) \to \sigma_0 \phi_t(s,0)$ - this also affects $u_{\nu}$ and $\partial_{\nu} u$ expansions.

\item (\textcolor{blue}{done}) Don't need to use Simon iteration argument since everything on $M$ now
\end{itemize}
\end{comment}
%
\section{Introduction} \label{IntroSection}
\begin{figure}[h!]
\centering
\includegraphics[scale=0.4]{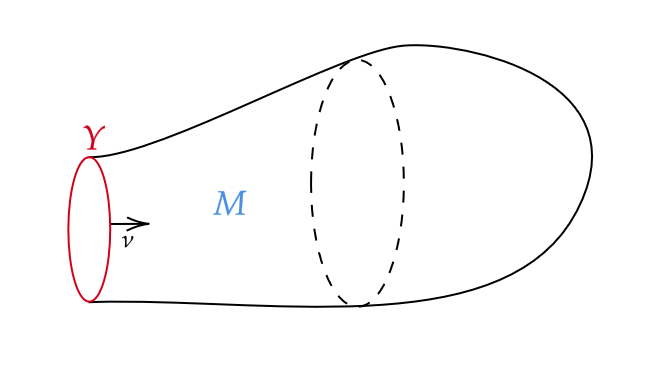}
\caption{Image of our set up}
\label{fig:acminimizationmainsetup}
\end{figure}

\begin{figure}[h!]
\centering
\includegraphics[scale=0.25]{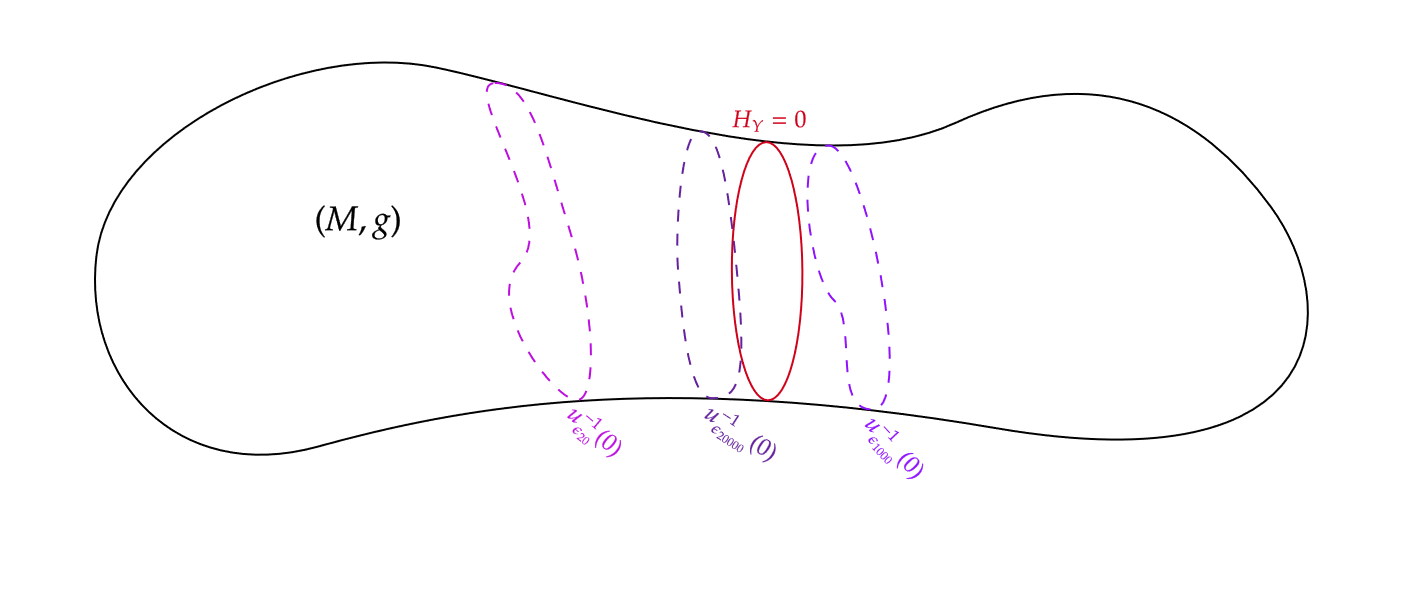}
\caption{Illustration of level set convergence to a minimal hypersurface}
\label{fig:aclevelsetconvergence}
\end{figure}

We work in $(M^n, g)$, a closed, smooth Riemannian manifold with boundary, $Y = \partial M$ (see \ref{fig:acminimizationmainsetup}). We assume that $Y$ is at least $C^{2,\alpha}$ and will state higher regularity when needed.
\nl
For any $\eps > 0$, there exists a nonnegative minimizer of the Allen--Cahn energy \cite{allen1979microscopic}
\begin{equation} \label{ACEnergy}
E_{\eps}(u) = \int_{M} \eps \frac{|\nabla_g u|^2}{2} + \frac{1}{\eps} W(u)
\end{equation}
such that $u\Big|_Y \equiv 0$. $W(u) = \frac{1}{4}(1 - u^2)^2$ can be take to be the standard double-well potential. Minimizers of this energy functional satisfy the Allen--Cahn equation on the interior of $M$
\begin{equation} \label{ACEquation}
\eps^2 \Delta_g u = W'(u) = u(u^2 - 1)
\end{equation}
On closed Riemannian manifolds, there is a well-known correspondence between zero sets of solutions to Allen--Cahn and minimal surfaces: Modica and Mortola (\cite{modica1985gradient} \cite{mortola1977esempio}) showed that the Allen--Cahn energy functional $\Gamma$-converges to perimeter. Under certain geometric constraints, Wang and Wei ([\cite{wang2019second}, Thm 1.1]) showed that the level sets of a sequence of stable solutions to \eqref{ACEquation}, $\{u_{\eps_i}\}$, converge to a minimal surface with good regularity. Similarly, given a minimal surface $Y \subseteq M$, Pacard and Ritore ([\cite{pacard2003constant}, theorem 4.1]) showed that one can construct solutions to Allen--Cahn with zero sets converging to $Y$ as $\eps \to 0$ (see \ref{fig:aclevelsetconvergence}).
\nl \nl
In this paper, we split figure \ref{fig:aclevelsetconvergence} and look at solutions on manifolds with boundary. We are concerned with the following questions: given the boundary $Y = \partial M$, how can one see the geometry of $Y$ in a solution, $u_{\eps, Y}$ to \eqref{ACEquation} with $0$ level set of $Y$? We take $u_{\eps, Y} = u_{\eps, Y}: M \to \R$ to be the non-negative minimizer of \eqref{ACEnergy} with zero Dirichlet data on $Y$. We then show an asymptotic expansion of the Neumann data, $\nu_Y(u_{\eps, Y})$, in powers of $\eps$, with the coefficients depending on the curvatures of $Y$. Finally, we apply our results to the setting of $(M^n, g)$ closed with $Y^{n-1} \subseteq M^n$ a separating hypersurface so that $M = M^+ \sqcup_{Y} M^-$.

\subsection{Background}
\begin{figure}[h!]
\centering
\includegraphics[scale=0.4]{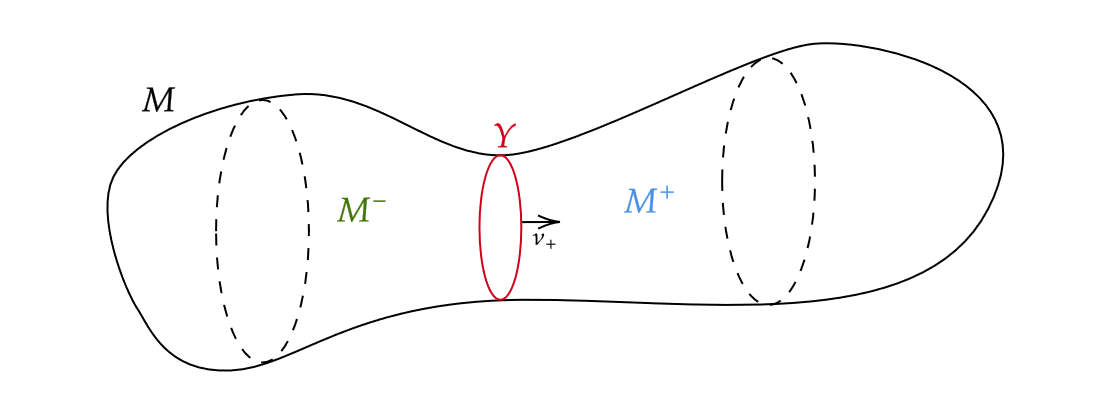}
\caption{Closed Setting}
\label{fig:acminimizationfullsetup}
\end{figure}

For $Y$ as above, consider the non-negative energy minimizer of \eqref{ACEnergy}, $u_{\eps, Y}$, with Dirichlet conditions on $Y$. By standard calculus of variations, this Dirichlet-minimizers exist. By work of Brezis--Oswald ([\cite{brezis1986remarks}, Thm 1]), there is at most one such solution to \eqref{ACEquation} on $M$ with this Dirichlet condition. Moreover, such a solution minimizes \eqref{ACEnergy} among all such functions. We ask, what is $\partial_{\nu} u_{\eps, Y}$? We describe this Neumann data, as well as an expansion of $u_{\eps, Y}$ itself \textbf{asymptotically in $\epsilon$}, by mimicking the techniques of Wang--Wei \cite{wang2019second} and also Mantoulidis [\cite{mantoulidis2022variational}, \S 4]. \nl \nl
Our main application is the closed setting of this problem. Let $(M^n, g)$ a smooth Riemannian manifold, and $Y^{n-1} \subseteq M^n$ a separating, two-sided hypersurface (see figure \ref{fig:acminimizationfullsetup}) 
\nl \noindent such that $M = M^+ \sqcup_Y M^-$. One can consider non-negative (resp. non-positive) minimizers of \eqref{ACEnergy} on $M^{\pm}$, referred to as $u_{\eps, Y}^{\pm}$. For $\nu = \nu_+$ a normal pointing inward to $M^+$, the condition $\partial_{\nu} u_{\eps, Y}^+ = \partial_{\nu} u_{\eps, Y}^-$ means that $u_{\eps, Y}^{\pm}$  can be pasted together to form a smooth solution to Allen--Cahn with level set on $Y$.  In particular, we are motivated by the following theorem of Pacard and Ritore [\cite{pacard2003constant}, Thm 1.1]: 
\begin{restatable}[Pacard--Ritore, Theorem 1.1]{thmm}{PacardRitoreTheorem} \label{PacardRitoreTheorem}
Assume that $(M, g)$ is an $n$-dimensional closed Riemannian manifold and $Y^{n-1} \subset M$ is a two sided, nondegenerate minimal hypersurface. Then there exists $\eps_0 > 0$ such that $\forall \eps \in (0, \eps_0)$ there exists $u_{\eps}$ solutions to the Allen--Cahn equation such that $u_{\eps}$ converges to $+1$ (resp. $-1$) on compact subsets of $(M^+)^o$ (resp. $(M^-)^o$) and 
\[
\mathscr{E}_{\eps}(u_{\eps}) \xrightarrow{\eps \to 0} \frac{1}{\sqrt{2}} A(Y)
\]
where $A(Y)$ is the $n-1$ dimensional area of $Y$.
\end{restatable}
\noindent We'll prove a theorem about the projection of the solutions constructed by Pacard and Ritore [\cite{pacard2003constant}, Thm 4.1] onto a specific kernel. \nl \nl
\noindent While Pacard and Ritore showed that one can construct solutions with zero sets converging to $Y$ minimal, Wang and Wei consider a stable sequence of solutions to \eqref{ACEquation}, $\{u_{\eps}\}$, and produce curvature bounds on the level sets [\cite{wang2019second}, Thm 1.1]. Recall that $B(u_{\eps})$ denotes the extended second fundamental form on graphical functions (see \cite{wang2019second}, Eq 1.4)
\begin{restatable}[Wang--Wei, Theorem 1.1]{thmm}{WWTheorem} \label{WWTheorem}
For any $\theta \in (0,1)$, $0 < b_1 \leq b_2 < 1$, and $\Lambda > 0$, there exist two constants $C = C(\theta, b_1, b_2, \Lambda)$ and $\eps_* = \eps(\theta, b_1, b_2, \Lambda)$ so that the following holds: suppose $u_{\eps}$ a stable solution of Allen--Cahn in $B_1(0) \subseteq \R^n$ satisfying 
\[
|\nabla u_{\eps}| \neq 0 \quad \text{and} \quad |B(u_{\eps})| \leq \Lambda \quad \text{in} \quad \{ |u_{\eps}| \leq 1 - b_2\} \cap B_1(0)
\]
If $n \leq 10$ and $\eps \leq \eps_*$, then for any $t \in [-1 + b_1, 1 - b_1]$, $\{u_{\eps} = t\}$ are smooth hypersurfaces and 
\[
[H(u_{\eps})]_{\theta} \leq C, \qquad ||H(u_{\eps})||_{C^0} \leq C \eps(\log |\log \eps|)^2
\]
where $H(u_{\eps})$ denotes the mean curvature of $\{u_{\eps} = t\}$
\end{restatable}
\noindent we are heavily inspired by the techniques used in their paper, though our results take on a different theme.
\subsection{Motivating Example }
\begin{figure}[h!]
\centering
\includegraphics[scale=0.3]{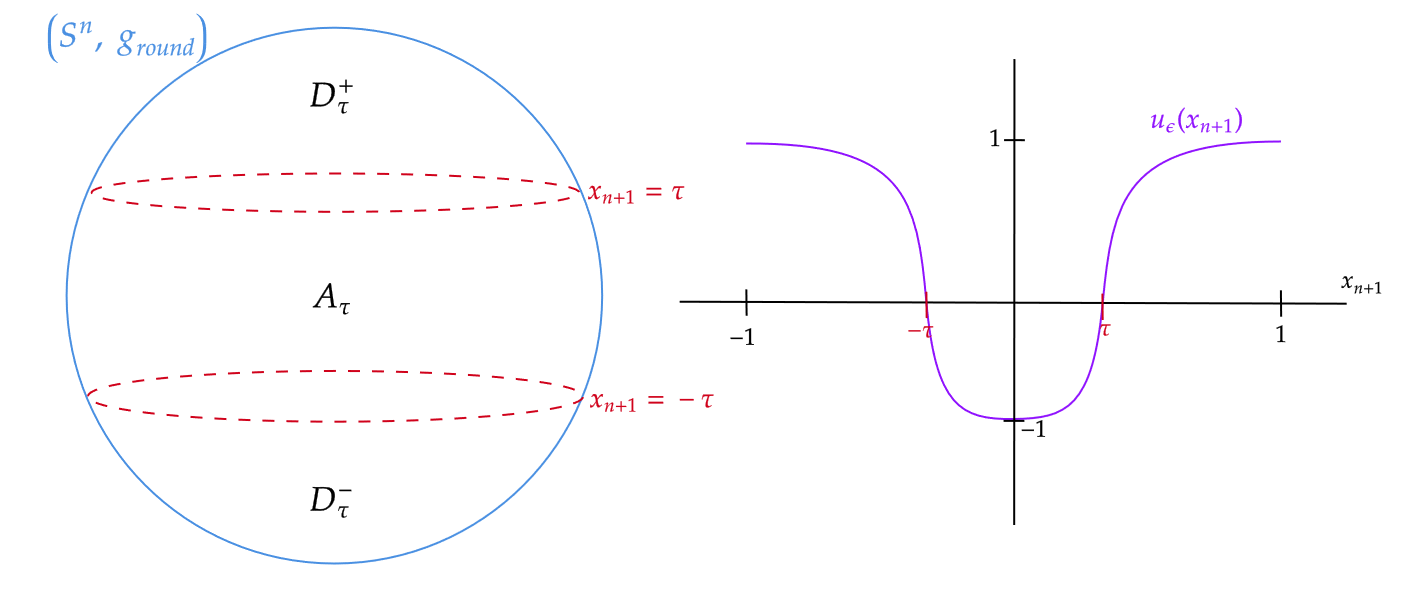}
\caption{Example of matching Neumann data on the sphere}
\label{fig:acsolutionsphere}
\end{figure}

As mentioned in the previous section, one can construct solutions of \eqref{ACEquation} by matching Dirichlet and Neumann conditions along a hypersurface. We're motivated by the following example from [\cite{guaraco2019min}, Ex. 19]: \nl \nl
\textbf{Example}: \; Let $M = S^n \subseteq \R^{n+1}$. Define regions $A_{\tau} = S^n \cap \{|x_{n+1}| < \tau\}$ and $S^n \backslash A_{\tau} = D_{\tau}^+ \cup D_{\tau}^-$, where $D_{\tau}^{\pm}$ are the discs forming the complement of the annulus $A_{\tau}$. Consider $u_{\eps, \tau}^{\pm}$ the \textit{nonnegative} energy minimizers of the Allen--Cahn energy on $D^{\pm}_{\tau}$. Let $v_{\eps, \tau}$ denote the \textit{nonpositive} energy minimizer on $A_{\tau}$ and define 
\[
\tilde{u}_{\tau, \eps}(p):= \begin{cases}
u_{\eps, \tau}^+(p) & p \in D_{\tau}^+ \\
u_{\eps, \tau}^-(p) & p \in D_{\tau}^- \\
v_{\eps, \tau}(p) & p \in A_{\tau} 
\end{cases}
\]
see figure \ref{fig:acsolutionsphere}. This is $C^0$ and a solution to Allen--Cahn on $S^n \backslash \partial A_{\tau}$. We aim to find $0 < \tau < 1$ such that $\tilde{u}_{\tau}$ is $C^1$ across $\partial A_{\tau}$, i.e. the Neumann data matches on $x_{n+1} = \pm \tau$. 
\nl \nl
\textbf{Proof Sketch}: One can show that
\[
C_{\eps, \tau}^{\pm} := \frac{\partial u_{\eps, \tau}^{\pm}}{\partial x^{n+1}} - \frac{\partial v_{\eps, \tau}}{\partial x^{n+1}} \Big|_{x_{n+1} = \pm \tau}
\] 
varies continuously with $\tau$ and is only dependent on $\eps$ and $\tau$ (i.e. these solutions are one dimensional). Note by symmetry that $C_{\eps, \tau}^- = -C_{\eps, \tau}^+$. \nl \nl
In particular for $\eps$ fixed and $\tau$ sufficiently close to $1$, $u_{\eps, \tau}^{\pm} \equiv 0$ while $v_{\eps, \tau} < 0$ on $\text{Int}(A_{\tau})$. Similarly for $\tau$ sufficiently close to $0$, $u_{\eps, \tau}^{\pm} > 0$ and $v_{\eps, \tau} \equiv 0$, i.e. for some $\delta > 0$
\begin{align*}
|1 - \tau| &< \delta \implies C_{\eps,\tau}^+ > 0 \\
|1 - \tau| &> 1 - \delta \implies C_{\eps,\tau}^+ < 0
\end{align*}
By continuity of $C_{\eps, \tau}^+$ there exists $\tau$ such that $C_{\eps, \tau}^{\pm} = 0$, and so that $\tilde{u}_{\tau, \eps}$ is a $C^1$ and hence smooth solution to Allen--Cahn on $S^n$. 
\subsection{Results} \label{Results}
\begin{figure}[h!]
\centering
\includegraphics[scale=0.4]{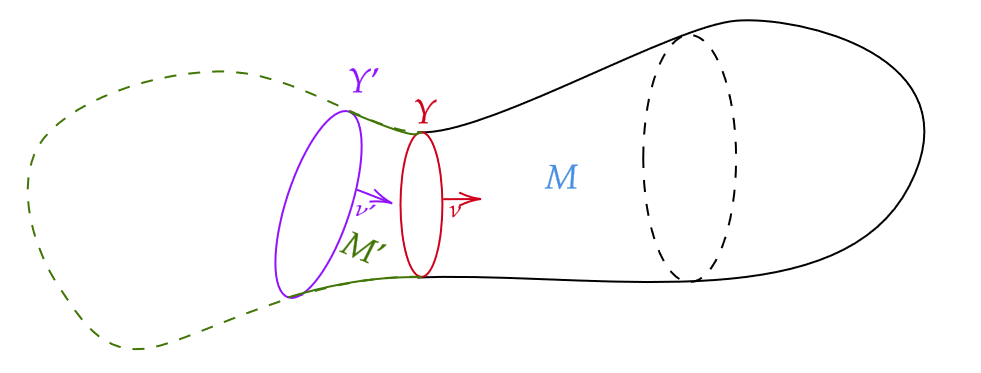}
\caption{Variable boundary visualization}
\label{fig:variableboundary}
\end{figure}

\noindent We recall the classic Dirichlet-to-Neumann map: consider an elliptic operator, $L$, arising from an energy functional. Given $(M, g)$ a Riemannian manifold and $\partial M$ smooth, consider $f : \partial M \to \R$. Suppose there exists a unique $\tilde{u}$ such that
\begin{align*}
\tilde{u} &: M \to \R \\
L \tilde{u} &= 0 \\
\tilde{u} \Big|_{\partial \Omega} &= f
\end{align*}
Then one can formulate a map 
\begin{align} \label{ClassicDToNMap}
\mathcal{D}&: H^{1}(\partial M) \to L^2(\partial M) \\ \nonumber
\mathcal{D}(f) &= \partial_{\nu}(\tilde{u})
\end{align}
see [\cite{mitrea1999boundary}, \S 7] for details. \nl \nl
Instead, we investigate the following Dirichlet-to-Neumann type map, where the Dirichlet data is fixed at $0$, but the manifold and its boundary, $(M, Y)$ are variable. This is pictured in figure \ref{fig:variableboundary}
where we think of $M$ as a subset of a larger closed manifold. Suppose we take the unique energy minimizer on $M$ such that
\begin{align*}
\tilde{u} &: M \to \R \\
\eps^2\Delta_g(\tilde{u}) - W'(\tilde{u}) &= 0 \\
\tilde{u}\Big|_{M} & > 0 \\
\tilde{u} \Big|_{Y} &= 0 
\end{align*}
Let $\nu$ be the positive normal to $Y$ and define
\begin{align*}
\mathcal{N}(Y) &:= \partial_{\nu}(\tilde{u}) \Big|_{Y}
\end{align*}
$\mathcal{N}(Y)$ is a global term and depends on the geometry of $M$ in a neighborhood of $Y$, as opposed to just on $Y$. We prove the following result which expands the Neumann data asymptotically in $\eps$:
\begin{restatable}{thmm}{NeumannResult}
\label{NeumannResult}
%Use $C^{3,\alpha}$ to get error term and show that $\Delta_t \phi$ is small
For $Y = \partial M$ a $C^{3,\alpha}$ hypersurface, and $u_{\eps, Y}$ the positive minimizer of $E_{\eps}$ on $M$ with $0$ Dirichlet condition on $Y$, we have that 
\[
\partial_{\nu} u_{\eps, Y} = \frac{1}{\eps \sqrt{2}} -  \frac{2}{3} H_Y + O(\eps^{1-\alpha})
\]
with error in $C^{\alpha}(Y)$.
\end{restatable}
\noindent \textbf{Remarks}:
\begin{itemize}
\item We can posit this more formally as a Dirichlet-to-Neumann operator when $M$ is closed and $Y \subseteq M$ is a separating hypersurface. Let $\mathcal{M}^{2,\alpha}(M)$ denote the space of $C^{2,\alpha}$ two sided, closed hypersurfaces with bounded geometry (cf. equation \eqref{UniformGeometry}). For $Y \subseteq \mathcal{M}^{2,\alpha}(M)$, let $U(Y)$ denote a normal neighborhood such that each $Y' \in U$ can be represented as 
\[
Y' = Y_{\eta} = F(Y, \eta) := \exp_Y(\eta(p) \nu(p))
\]
with $\eta \in C^{2,\alpha}(Y)$ and $\nu$ a normal to $Y$. Define
\begin{align} \label{FMap}
F&: Y \times (-\delta, \delta) \to M \\ \nonumber
F(p,t) &= \exp_p (t \nu(p)) \\ \nonumber
F_{\eta}(p, t)& = \exp_p((t + \eta(p)) \nu(p))
\end{align}
so that
\begin{align} \label{DtoNMap}
\mathcal{N} &: U(Y) \to C^{1,\alpha}(Y) \\ \nonumber
Y_{\eta} & \mapsto F_{\eta}^*(\nu_{Y_{\eta}} u_{\eps, Y_{\eta}}^+ ) \Big|_{t = 0} \\ \nonumber
& = (F_{\eta}^{-1})_*(\nu_{\eta}) F_{\eta}^*(u_{\eps, Y_{\eta}}^+) \Big|_{t = 0}
\end{align}
Now noting that $U(Y) \cong \{\eta\} =: V \subseteq C^{2,\alpha}(Y)$, we can frame $\mathcal{N}: V \to C^{\alpha}(Y)$ as a map between functions on $Y$. In this sense, the variable initial level set, $Y_{\eta}$, is the Dirichlet data, and the normal derivative is the Neumann data, which forms a Dirichlet-to-Neumman type map on $Y$ after pulling back. Our results do not just apply to perturbations of a fixed $Y$, as all of our theorems (with the exception of \cref{ProjectionTheorem}) apply to \textit{any} separating hypersurface of the appropriate regularity and bounded geometry. For this, we use the term ``Dirichlet-to-Neumann \textit{type} map" to describe $\mathcal{N}$
\item As noted in [\cite{guaraco2019min}, Ex. 8], for $Y$ fixed, a positive minimizer on $M^{\pm}$ will exist for $\eps$ sufficiently small.
\end{itemize}
\noindent We begin with a key estimate on $L_{\eps} = \eps^2 \Delta_g - W''(\bg_{\eps})$:
\begin{restatable}{thmm}{SchauderBoundLemma} \label{SchauderBoundLemma}
Let $Y = \partial M^+$ a $C^{2,\alpha}$ surface, and suppose $f: M^+ \to \R$ in $C^{2,\alpha}_{\eps}(M^+)$ with $f(s,0) \equiv 0$. Then there exists an $\eps_0 > 0$ sufficiently small, independent of $f$, such that for all $\eps < \eps_0$, we have 
\[
||f||_{C^{2,\alpha}_{\eps}(M^+)} \leq K ||L_{\eps} f||_{C^{\alpha}_{\eps}(M)}
\]
for $K$ independent of $\eps$.
\end{restatable}
\noindent Immediately, we see that the linearized Allen-Cahn operator is invertible as a map from $C^{2,\alpha}_{\eps}(M^+) \cap C_0(M^+)$ (i.e. zero boundary conditions) to $C^{\alpha}_{\eps}(M^+)$ (see \cite{gilbarg1977elliptic}, Theorem 6.15). \nl \nl
\vspace{0.2cm}\noindent After establishing this, theorem \ref{NeumannResult} is then proved in the following manner:
\begin{enumerate}
\item Let $t$ denote the signed distance from $Y$ and $s$ a fermi coordinate on $Y$. We decompose 
\[
u_{\eps, Y}(s,t) = \bg(t/\eps) + \phi(s,t)
\]
where $\bg(t)$ is a modification of the heteroclinic solution, and we rewrite the Allen--Cahn equation in terms of $\phi$
\item We prove a modified Schauder estimate, 
\[
||\phi||_{C^{2,\alpha}_{\eps}(M^+)} \leq ||L \phi||_{C^{\alpha}_{\eps}(M^+)}
\]
reminiscent of [\cite{pacard2012role}, Prop $3.21$]

\item We integrate the Allen--Cahn equation, showing that 
\[
\partial_t\phi(s,0)\Big|_{t = 0} = H_Y(s) \sigma_0 + \int_0^{-\omega\eps \ln(\eps)} (\Delta_t \phi) \dot{\bg}_{\eps} + O(\eps^2)
\]

\item We show that $\int_0^{-\omega\eps \ln(\eps)} (\Delta_t \phi) \dot{\bg}_{\eps}$ is small by proving better $C^{2,\alpha}_{\eps}$ estimates for $\partial_{s_i} \phi$. This mimics [\cite{wang2019second}, \S 7]

\end{enumerate}
\noindent We can improve our analysis of the Neumann data when $H_Y = 0$:
\begin{restatable}{thmm}{ImprovedNeumannMinimal} \label{ImprovedNeumannMinimal}
When $Y = \partial M$ is minimal, we have that 
\[
\nu^+(u_{\eps, Y}^+) = \frac{1}{\eps \sqrt{2}} + \sigma_0^{-1} \kappa_0 \eps [\Ric_Y(\nu, \nu) + |A_Y|^2] + O(\eps^{2-\alpha})
\]
with error in $C^{\alpha}(Y)$.
\end{restatable}
\noindent In the manifold with boundary setting, we see that as $Y$ is more regular, we can capture more terms in the expansion of $\nu(u_{\eps, Y})$. This culminates in our main theorem
\begin{restatable}{thmm}{InductiveExpansion} \label{InductiveExpansion}
For $Y = \partial M$ a $C^{\overline{k} + 3,\alpha}$ hypersurface, the minimizer of \eqref{ACEnergy} can be expanded as 
\begin{align} \label{InductiveExpansionEquation}
u_{\eps}(s,t) &= \bg_{\eps}(t) + \sum_{i = 1}^{\overline{k}} \eps^i \cdot \left( \sum_{j = 0}^{M_i} a_{i,j}\left(\{\partial_s^{\beta} \partial_t^j H_t(s) |_{t = 0}\}_{j + |\beta| \leq i}\right)\overline{w}_{i,j,\eps}(t) \right) + \phi \\ \nonumber
||\phi||_{C^{2,\alpha}_{\eps}(M)} & = O( \eps^{\overline{k} + 1})
\end{align}
where $a_{i,j}(s)$ are a collection of polynomials in derivatives of $H_t$ up to a certain order and $\overline{w}_{i,j,\eps}(t) = \overline{w}_{i,j}(t/\eps)$ are modifications of functions $w_{i,j,\eps}(t) = w_{i,j}(t/\eps)$ satisfying
\begin{align*}
w_{i,j}&: [0, + \infty) \to \R \\
||w_{i,j}||_{C^{\infty}_{loc}} &= O(1)
\end{align*}
and are exponentially decaying in $C^{\infty}$. 
\end{restatable}
\noindent When $\overline{k} \geq 1$ this yields
\begin{restatable}{corr}{NextOrderCorollary} \label{NextOrderCorollary}
For $u_{\eps}$ a solution to Allen--Cahn with Dirichlet data on $Y$ a $C^{4,\alpha}$ hypersurface, we have that 
\begin{align*}
u_{\eps, Y}^+(s,t) &= \bg_{\eps}(t) + \eps H_Y(s) \bw_{\eps}(t) + \phi(s,t) \\
||\phi||_{C^{2,\alpha}_{\eps}(M)} &= O(\eps^2)
\end{align*}
\end{restatable}
\noindent When $Y$ is minimal (and has no singular set by assumption at beginning of \S \ref{IntroSection})
\begin{restatable}{corr}{NextOrderMinimalCorollary} \label{NextOrderMinimalCorollary}
For $u_{\eps, Y}$ a solution to Allen--Cahn with Dirichlet data on $Y$, a minimal surface, the expansion in \eqref{InductiveExpansionEquation} exists to any order.
\end{restatable}
\noindent \rmk \; In general, we see that both the expansion of $u_{\eps, Y}$ and its neumann derivative are \textit{asymptotically local} in terms of a series expansion in powers of $\eps$, despite being global quantities determined by the geometry of all of $M$, not just the geometry in a neighborhood of $Y$. \nl \nl
\noindent \rmk \; Similar expansions have been done by Wang--Wei \cite{wang2019finite}, Chodosh--Mantoulidis \cite{chodosh2020minimal}, and Mantoulidis \cite{mantoulidis2022variational} among other authors. These works begin with $u: M \to \R$ a smooth solution to \eqref{ACEquation} and then expand $u$ about its zero-set. This approach actually gives the zero-set better regularity by a Simons-type equation (see [\cite{wang2019finite}, Lem 8.6]), allowing for more terms in the expansion. By contrast, we start with $Y$, a prescribed zero set with limited regularity, and $u_{\eps, Y}^{\pm}: M^{\pm} \to \R$ one sided solutions, for which the Simons-type equation does not apply. \nl \nl
\noindent Returning to the setting of $M$ closed and $Y \subseteq M$ a separating hypersurface: take $\eta \in C^{2,\alpha} (Y)$, along with \eqref{FMap}
\begin{align*}
F_{\eta}&: Y \to U(Y) \subseteq M \\
F_{\eta}(p,0) &:= \exp_p(\eta(p) \nu(p)) \\
Y_{\eta} &= \{F_{\eta}(p) \; | \; p \in Y\}
\end{align*}
\begin{figure}[h!]
\centering
\includegraphics[scale=0.5]{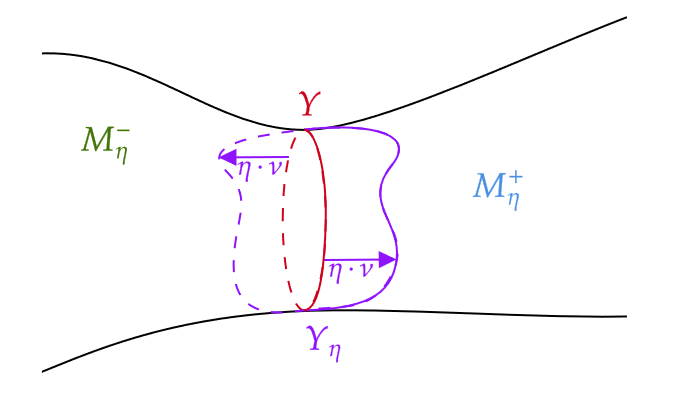}
\caption{Perturbation of $Y$ and the corresponding splitting of $M = M_{\eta}^+ \sqcup_{Y_{\eta}} M_{\eta}^-$}
\label{fig:acetaperturbation}
\end{figure}
\noindent where $U(Y)$ is some open neighborhood of $Y$ in $M$. Decompose $M = M_{\eta}^+ \sqcup_{Y_{\eta}} M_{\eta}^-$ (see figure \ref{fig:acetaperturbation}) and consider the (positive) energy minimizers $u_{\eps, \eta}^{\pm}$ on $M_{\eta}^{\pm}$. By Brezis--Oswald [\cite{brezis1986remarks}, Thm 1], these are the unique solutions to \eqref{ACEquation} on $M_{\eta}^{\pm}$ and we can paste them together to form:
\begin{align} \label{PastedSolution}
u_{\eps, \eta}: &M \to \R \\ \nonumber
u_{\eps, \eta} &= \begin{cases}
u_{\eps, \eta}^+(p) & p \in M_{\eta}^+ \\
-u_{\eps, \eta}^-(p) & p \in M_{\eta}^-
\end{cases}
\end{align}
We now use the map in \eqref{DtoNMap} 
\begin{align*}
\mathcal{N}^{\pm}&: U(Y) \subseteq C^{2,\alpha}(Y) \to C^{1,\alpha}(Y) \\
\mathcal{N}^{\pm}(\eta) &= F_{\eta}^*(\partial_{\nu^{\pm}_{Y_{\eta}}} u_{\eps, \eta}^{\pm})(p)
\end{align*}
Note that the Neumann data of $u_{\eps, \eta}$ matching along $Y_{\eta}$ is equivalent to 
\[
\mathcal{N}^+(\eta) - \mathcal{N}^{-}(\eta) = 0
\]
When this is the case, $u_{\eps, \eta}$ is a smooth solution to \eqref{ACEquation} and we can characterize the projection of $u_{\eps, \eta}$ onto $\dot{\bg}_{\eps}$, the kernel of $L_{\eps} := \eps^2 \partial_t^2 - W''(g_{\eps})$:
\begin{restatable}{thmm}{ProjectionTheorem} \label{ProjectionTheorem}
Let $Y \subseteq M$ a minimal separating hypersurface in a closed, smooth Riemannian manifold. For $u_{\eps, \eta}:M \to \R$ as in \eqref{PastedSolution}, suppose that $u_{\eps, \eta}$ is $C^1$ across $Y_{\eta}$ and $||\eta||_{C^{2,\alpha}} = O(\eps^{1 + \beta})$ for some $\beta \geq \alpha$ fixed. Then
\[
\int_{\R} \Delta_Y(\phi) \dot{\bg}_{\eps}(t) dt = 2 \sigma_0 J_Y(\eta) + \tilde{O}(\eps^{1 + 2 \beta})
\]
with error holding in $C^{\alpha}(Y)$. If we further have that $Y$ is non-degenerate, then
\[
\int_{\omega\eps \ln(\eps)}^{-\omega\eps \ln(\eps)} u_{\eps, \eta}(s,t) \dot{\bg}_{\eps}(t) dt = \frac{\sqrt{2}}{3} \eta(s) + \tilde{O}(\eps^2, \eta^2)
\]
with error in $C^{2,\alpha}(Y)$
\end{restatable}
\noindent \rmk \; The above theorem tells us that when we perturb $Y \to Y_{\eta}$ to find a solution to \eqref{ACEquation} with zero set $Y_{\eta}$, then we can detect $\eta$ via the projection of our solution onto $\dot{\bg}_{\eps}$. Also note that for $u_{\eps, \eta}(s,t) = \bg_{\eps}(t) + \phi(s,t)$, 
\[
\int_{\omega\eps \ln(\eps)}^{-\omega\eps \ln(\eps)} u_{\eps, \eta}(s,t) \dot{\bg}_{\eps}(t) dt = \int_{\omega\eps \ln(\eps)}^{-\omega\eps \ln(\eps)} \phi(s,t) \dot{\bg}_{\eps}(t) dt + O(\eps^k)
\]
so \cref{ProjectionTheorem} is equivalent to computing the projection of $\phi$ onto $\dot{g}_{\eps}$. Further note that \cref{ProjectionTheorem} differs from corollary \ref{NextOrderCorollary} in that we compute a two-sided integral for theorem \ref{ProjectionTheorem}.

\subsection{Paper Organization}
This paper is organized as follows:
\begin{itemize}
\item In \S \ref{SetUp}, we define notation and recall some known geometric equations and quantities.

\item In \S \ref{NormalDerivativeSection}, we pose the initial decomposition of $u_{\eps, Y}^+ = \bg(t/\eps) + \phi$ where $\bg$ is a modification of the heteroclinic. We prove \cref{SchauderBoundLemma}, which is used to estimate $||\phi||_{C^{2,\alpha}_{\eps}}$ and $||\partial_{s_i} \phi||_{C^{2,\alpha}_{\eps}}$. We then prove \cref{NeumannResult}

\item In \S \ref{HigherOrderSection}, we prove higher order expansions of the Neumann data, \cref{ImprovedNeumannMinimal} and our main result, \cref{InductiveExpansion}. This theorem says that given $Y$ a $C^{k+3,\alpha}$ surface with bounded geometry, we can expand $u_{\eps, Y}^{\pm}$ to order $\eps^{k}$.

\item Finally, in \S \ref{ProofOfProjectionSection}, we prove \cref{ProjectionTheorem}
\end{itemize}
\subsection{Acknowledgements and Dedication}
The author would like to thank Otis Chodosh for presenting him with this problem and many fruitful conversations over the course of several months. The author would also like to thank Rafe Mazzeo and \'{E}rico Silva for their mathematical perspectives on the project. Furthermore, the author would like to thank Yujie Wu and Shuli Chen for their feedback and companionship. \nl \nl
\noindent The author dedicates this project to his grandmother, Shirley Kuo, who endured many hardships to come to the US for a better life. In the US, she was denied an opportunity to do a PhD due to sexism, despite being overqualified. This paper is in honor of her.
\section{Set up} \label{SetUp}
We first describe the manifold with boundary setting. Let $(M^{n}, g)$ a riemannian manifold with $Y = \partial M$ a $C^{k,\alpha}$ surface for some $k \geq 2$. Throughout this paper, we'll assume uniform bounds on the geometry of $Y$, i.e. we fix a $C > 0$ (independent of $\eps$) such that
\begin{equation} \label{UniformGeometry}
||A_Y||_{C^{k-2,\alpha}} \leq C
\end{equation}
We define $u_{\eps, Y}$ (also notated as $u_{\eps}$) to be the energy minimizer of \eqref{ACEnergy} on $M$ with Dirichlet conditions on $Y$. Let $p \in Y$ a base point, $\{E_i\}$ an orthonormal frame for $TY$ at $p$, and $\nu$ a normal on $Y$ with respect to $g$. We coordinatize a tubular neighborhood of $Y$ via the following maps
\begin{align} \label{CoordinatizationMap}
G&: B_1(0)^{n-1} \to Y \\ \nonumber
G(s) &= \exp_{p}^Y(s^i E_i) \\ \label{FMapGeneral}
F&: B_1(0)^{n-1} \times [0, \delta_0) \to M \\ \nonumber
F(s,t)&:= \exp_{G(s)}^{M}(\nu(G(s)) t)
\end{align}
for any $\omega > 5$ fixed and finite. Here, $\exp^Y_p: B_1(0) \to Y$ denotes the exponential map into $Y$, and $\exp^{M}: Y \times [0, \delta_0) \to M$ is the $M$-exponential map. We can coordinatize $u = u(s, t)$ in an $\eps$-neighborhood of $Y$ via the above. In this neighborhood, we can expand the metric and second fundamental form, $g_{ij}(s,t)$ and $A_{ij}(s,t)$, in coordinates, smoothly in $t$ from the following equations
\begin{align} \label{metricExpansion}
g_{ij}(s,t) &= g_{ij}(s,0) - 2 t \sum_{k = 1}^n A_i^k(s,0) g_{jk}(s,0) + t^2 \sum_{k, \ell = 1}^n A^2(\partial_i, \partial_j) - \text{Rm}_g(\partial_i, \partial_t, \partial_t, \partial_j) + O(t^3) \\ \label{secondFundamentalExpansion}
A(s,t) &= A(s,0) + t \left[A^2|_{(s,r)} - \text{Rm}(\cdot, \partial_t, \partial_t, \cdot) |_{(s,r)}\right] + O(t^2) \\
H(s,t) &= H(s,0) + t [-|A|^2|_{(s,t)} - \Ric_g(\partial_t, \partial_t)|_{(s,t)}] + O(t^2)
\end{align}
Here, $g(s,0)$, $A(s,0)$, and $H(s,0)$ denote the corresponding geometric quantities on $Y$. Moreover, $A^2$ denotes a single trace of $A \otimes A$ (see \cite{chodosh2020minimal} A.1, A.2). We can decompose the laplacian on $M$ in a neighborhood of $Y$ via
\begin{equation} \label{LaplacianDecomposition}
\Delta_g = \Delta_t - H_t \partial_t + \partial_t^2
\end{equation}
where $\Delta_t = \Delta_{Y + t \nu}$ is the laplacian on the surface $Y + t \nu = \{p \; | \; d_{\text{signed}}(p, Y) = t\}$ and $H_t$ denotes the mean curvature of $Y + t \nu$. See [\cite{wang2019second}, \S $3$] for details. In light of this notation, $H_t \Big|_{t = 0} = H_0 = H_Y$ and we'll use $H_0$ and $H_Y$ interchangeably. Similarly, $\Delta_t \Big|_{t = 0} = \Delta_0 = \Delta_Y$ and we'll use these two interchangeably as well. While the above expansions hold on $t \in [0, \delta_0)$, we will often restrict $u$ to $t \in [0, -\omega\eps \ln(\eps))$ as the behavior of $u$ is well-understood for $t > -\omega\eps \ln(\eps)$. In the closed setting for which $Y \subseteq M$ is separating, we decompose $M = M^+ \sqcup_Y M^-$ and use the above framework for $(M^+, Y = \partial M^+)$ and $(M^-, Y = \partial M^-)$ respectively. \nl \nl
We also define
\begin{equation} \label{rescaledMetric}
g_{\eps}:= \eps^{-2} g
\end{equation}
along with the geometric H\"older spaces with respect to $g_{\eps}$, i.e. 
\begin{align*}
||f||_{C^{0, \alpha}_{\eps}}&:= ||f||_{C^0} + [f]_{0,\alpha, \eps} \\
[f]_{\alpha, M} &= \sup_{p_1 \neq p_2 \in M } \frac{|f(p_1) - f(p_2)|}{|\text{dist}_{g}(p_1, p_2)|^{\alpha}} \\ 
[f]_{k, \alpha, \eps} &:= [f]_{k, \alpha, \eps, M} = \sup_{\beta \st |\beta| = k} \; \sup_{p_1 \neq p_2 \in M} \frac{|D^{\beta} f(p_1) - D^{\beta} f(p_2)|}{\text{dist}_{g_{\eps}}(p_1, p_2)^{\alpha} } \\
||f||_{C^{k, \alpha}_{\eps}}&:= \sum_{j = 0}^k ||D^j f||_{C^{0, \alpha}_{\eps}} \\
C^{k,\alpha}_{\eps, 0}(M) &= C^{k,\alpha}_{\eps} \bigcap \Big\{f: M \to \R \; | \; f\Big|_{\partial M} \equiv 0 \Big\}
\end{align*}
Note that equation \eqref{ACEquation} becomes
\[
\Delta_{g_{\eps}} u = W'(u)
\]
i.e. by rescaling the metric, we can set $\eps = 1$. In accordance with this, we can define the following blow up maps:
\begin{align} \label{BlowUpMaps}
G_{\eps}&: B_{\eps^{-1}}(0)^{n-1} \to Y \\ \nonumber
G_{\eps}(\sigma) &:= \exp_{p, g_{\eps}}^Y(\sigma^i E_i) \\ \nonumber
F_{\eps}&: B_{\eps^{-1}}(0)^{n-1} \times [0, - \omega\ln(\eps)) \to M \\ \nonumber
F_{\eps}(\sigma,\tau)&:= \exp_{G(\sigma), g_{\eps}}^{NY}(\nu(G(\sigma)) \tau)
\end{align}
for any $\omega > 5$ fixed and finite. We may refer to $(\sigma, \tau)$ as ``scaled" fermi coordinates, as opposed to $(s,t)$, the actual fermi coordinates. \nl \nl 
For $\eta \in C^{2,\alpha}(Y)$ non-negative, define the perturbed graph
\begin{equation} \label{EtaPerturbationEquation}
Y_{\eta} := \{F(s, \eta(s)) \; | \; s \in Y\}
\end{equation}
where $F(s,\eta(s))$ (see \eqref{FMapGeneral}) is the unique point a signed distance of $\eta(s)$ away from $G(s) \in Y$. As with $Y$, $Y_{\eta} = \partial M_{\eta}$, where $M_{\eta}$ consists of all points a non-negative distance from $Y_{\eta}$. We then define $u_{\eta, \eps}$, the minimizer of $E_{\eps}$ on $M_{\eta}$ with $0$ Dirichlet condition on $Y_{\eta}$. \nl \nl 
For the closed setting, we have $(M^n, g)$ a closed riemannian manifold, and $Y^{n-1} \subseteq M^n$, a separating, two-sided hypersurface.  In this case $\eta$ in \eqref{EtaPerturbationEquation} is real-valued. Moreover, $Y_{\eta}$ divides $M$ into $M^+_{\eta}$ and $M^-_{\eta}$. We then define $u^{\pm}_{\eta, \eps}$, the non-negative (resp. non-positive minimizers) of $E_{\eps}$ on $M_{\eta}^{\pm}$ with $0$ Dirichlet condition on $Y_{\eta} = \partial M^{\pm}_{\eta}$ (see figure \ref{fig:acetaperturbation})
\subsection{Constants and Definitions}
We recall the $1$-dimensional solution to equation \eqref{ACEquation}, the \textit{heteroclinic solution}, denoted by 
\[
g(t) := \tanh(t/\sqrt{2})
\]
This notation follows previous convention, but the author notes the abuse of notation between the heteroclinic solution and the metric. The context should make it clear when one is used over the other. \nl \nl
For any $\eps > 0$ we define 
\[
\bg(t) := [1 - \chi(t/(-\omega \ln(\eps)))] g(t) + \chi(t/(-\omega \ln(\eps)))
\]
where $\chi(t)$ is smooth function which is $0$ for $t < 1$, goes from $0 \to 1$ on $[1,2]$, and is $1$ for $t \geq 2$. For reference, we'll also use $\chi_{\delta}(t):= \chi(t/\delta)$. Note that 
\[
\partial_t^2 \bg - W'(\bg) = O(\eps^{\omega})
\]
in a $C^{k,\alpha}(\R^+)$ sense and is supported on $[-\omega \ln(\eps), -2 \omega \ln(\eps)]$. Now let $\dot{g}$, $\ddot{g}$ denote the first and second derivatives of $g(t)$. We further denote
\[
g_{\eps}(t):= g(t/\eps), \qquad \dot{g}_{\eps}:= \dot{g}(t/\eps), \qquad \ddot{g}_{\eps}(t):= \ddot{g}(t/\eps)
\]
to be the rescaled versions of $g$ and its derivatives and analogously for $\bg_{\eps}, \dot{\bg}_{\eps}, \ddot{\bg}_{\eps}$. Furthermore, let 
\[
\Rwe = \eps^2 \partial_t^2 \bg_{\eps} - W'(\bg_{\eps})
\]
which is $O(\eps^{\omega})$ in $C^{k,\alpha}_{\eps}(\R)$ and supported on $[- \omega \eps \ln(\eps), -2 \omega \eps \ln(\eps)]$.
\nl \nl
%
%\begin{figure}[h!]
%\centering
%\includegraphics[scale=0.65]{AC_Energy_Minimization_Pictures/Heteroclinic_and_derivatives}
%\caption{}
%\label{fig:heteroclinicandderivatives}
%\end{figure}
%
%\nl
\noindent Define the constants 
\[
\sigma_0:= \int_0^{\infty} \dot{\bg}^2 dt = \frac{\sqrt{2}}{3}, \qquad \kappa_0:= \int_0^{\infty} t \dot{\bg}^2 dt = \frac{1}{6} [4 \ln(2) - 1], \qquad \sigma = \dot{g}(0) = \frac{1}{\sqrt{2}}
\]
Similarly, consider $w: [0, \infty) \to \R$, the solution to
\begin{align} \label{wEquation}
w''(t)  - W''(g) w(t) &= g(t) \\
w(0) &= 0 \\
\lim_{t \to \infty} w(t) &= 0
\end{align}
which exists and is unique by section \S \ref{HalfLineSection} in the appendix (see also [\cite{alikakos1996critical}, Lemma B.1, remark B.3]). We note that $\dot{w}(0) < 0$. As with $\bg$ and $\bg_{\eps}$, let 
\[
\bw(t) := [1 - \chi(t/(-\omega \ln(\eps)))] \overline{w}(t)
\]
(i.e. smooth cut off to $0$). Also let $w_{\eps}(t):= w(t/\eps)$ and similar for $\dot{w}_{\eps}, \ddot{w}_{\eps}, \dot{\bw}_{\eps}, \ddot{\bw}_{\eps}$. In general, for any exponentially decaying function satisfying an ODE similar to equation \eqref{wEquation}, we'll adopt the same notation of 
\[
f \to \overline{f}(t) = (1 - \chi(t/(-\omega \ln(\eps)))) f(t)
\]
With this, we define the linearized Allen--Cahn operator about $\bg_{\eps}$:
\[
L_{\eps} := \eps^2\Delta_g - W''(\bg_{\eps}): C^{2,\alpha}(M^+) \to C^{\alpha}(M^+)
\]
We also define big $O$ and $\tilde{O}$ notation to capture the size of error terms. We write $E = O(\eps^m)$ or $||E||_{C^{k,\alpha}_{\eps}} = O(\eps^m)$ to denote
\[
||E||_{C^{k,\alpha}_{\eps}} \leq C \eps^m
\] 
for some $C$ independent of $\eps$. Similarly, $E = \tilde{O}(f_1, f_2, \dots)$ denotes error depending on functions:
\begin{equation} \label{TildeErrorNotation}
E = \tilde{O}(f_1, f_2, \dots) \implies ||E||_{C^{k,\alpha}_{\eps}} \leq C \sum_i ||f_i||_{C^{k,\alpha}_{\eps}}
\end{equation}
for some $C$ independent of $\eps$ and the $\{f_i\}$. \nl \nl
Finally, we establish a definition for exponentially decaying functions:
\begin{definition}
A function, $f: [0, \infty) \to \R$, is exponentially decaying if there exists $\gamma > 0$, $C > 0$, $t_0 > 0$ such that
\[
\forall t > t_0, \qquad |f(t)| \leq C e^{-\gamma t}
\]
Moreover, $f$ is exponentially decaying in $C^k$ if such a bound holds for all $k$ derivatives of $f$. $f$ is exponentially decaying in $C^{\infty}$ if there exists such an $\gamma_k$, $C_k$, $t_k$ such that the above holds for each derivative of $f$.
\end{definition}
\section{Normal Derivative for $Y$} \label{NormalDerivativeSection}
In this section we work in the manifold with boundary setting. The goal is to show \cref{NeumannResult}
\NeumannResult*
\subsection{Initial Decomposition} \label{InitialDecomposition}
Decompose
\begin{equation} \label{PhiDecomposition}
u_{\eps}(s,t) = \bg_{\eps}(t) + \phi(s,t)
\end{equation}
with the above holding on a tubular neighborhood of $Y$ in normal coordinates i.e. $Y \times [0, \delta_0)$. We recall the following initial bound that $||\phi||_{C^{k,\alpha}_{\eps}} = o(1)$ as $\eps \to 0$:
\begin{lemma} \label{InitialPhiBound}
Let $Y$ be a $C^{k+1,\alpha}$ surface. For $\phi(s,t): M^+ \to \R$ as in \eqref{PhiDecomposition}, we have that for all $\mu > 0$, there exists an $\eps_0(\mu)$ such that for all $\eps < \eps_0$,
\[
||\phi||_{C^{k,\alpha}_{\eps}(M^+)} \leq \mu
\]
\end{lemma}
\noindent \Pf Let $R > 0$ to be determined. Note that $\phi: B_1(0)^{n-1} \times [0, \eps R) \to \R$ is smooth away from the boundary and $C^{k+1,\alpha}$ near the boundary (i.e. about $t = 0$, see [\cite{gilbarg1977elliptic}, Lemma 6.18]). On this subdomain, we have
\begin{equation} \label{PhiDecay}
\phi = o(1) \in C^{k,\alpha}_{\eps}(Y \times [0, \eps R])
\end{equation}
as $\eps \to 0$. This follows by 
\begin{itemize}
\item Blowing up our sequence of $u_{\eps}$ on $(M, g_{\eps})$, to get a $C^{k+1}$ solution, $u_{\infty}: \R^{n-1} \times \R^+ \to \R$ to \eqref{ACEquation} with $\eps = 1$
\item Considering the odd reflection of $u_{\infty}$, to get $\tilde{u}_{\infty}: \R^n \to \R$ a smooth solution on the whole space (the dirichlet and neumann data match at $t = 0$!)
\item Using a classification of solutions to Allen--Cahn on $\R^{n}$ with $u^{-1}(0) = \{x_{n} = 0\}$ [\cite{hamel2021half}, Thm 3]
\end{itemize}
This gives us $C^{k+1}_{loc}$ convergence on $\R^{n} = \R^{n-1} \times \R$, which when we restrict to $t \in [0, R]$ gives uniform convergence in $t$ for any $B_r(p) \subseteq \R^{n-1}$ fixed. Bounding $C^{k,\alpha}$ norms by $C^{k+1}$ norms, translating this back to the unscaled setting of $Y \times [0, \eps R]$, and noting that $Y$ is closed, we see that \eqref{PhiDecay} holds. For $t > R \eps$, we recall that for $u_{\eps}$ a solution to \eqref{ACEquation} with $u_{\eps}^{-1}(0) = Y$, we have the following decay estimate from [\cite{guaraco2019min}, Exercise 10, Remark] for any $\ell \in \Z^+$, $\alpha >0$,
\begin{equation} \label{PhiDistanceDecay}
||\phi||_{C^{\ell,\alpha}} = ||u_{\eps}(s,t) - g(t/\eps)||_{C^{\ell,\alpha}} \leq C e^{- \sigma t/\eps}
\end{equation}
where $C = C(\ell,\alpha)$ and $\sigma > 0$ independent of $\eps, \; \ell, \; \alpha$. This holds for all $t$, and on $I_2 = (\eps R, \delta_0)$ for $R$ sufficiently large,
\[
||\phi||_{C^{k,\alpha}(\{t > R \eps\})} \leq C e^{- \sigma R} \leq \mu
\]
for some $\eps$ sufficiently small. Since $||f||_{C^{k,\alpha}_{\eps}} \leq ||f||_{C^{k,\alpha}}$, the conclusion follows. \qed \nl \nl
Having given an initial bound on $\phi$, we write a PDE describing it
\begin{lemma} 
For $\phi(s,t)$ as in \eqref{PhiDecomposition}, we have 
\begin{equation} \label{ACDecomposition}
L_{\eps}(\phi) = \eps H_t \dot{\bg}_{\eps} + Q_0(\phi)
\end{equation}
where $Q_0(\phi) = \tilde{O}(\phi^2, \eps^{\omega})$ holds in $C^{k,\alpha}_{\eps}$ for any $(k, \alpha)$
\end{lemma}
\noindent \Pf
In this decomposition, the Allen--Cahn equation for $t \in [0, \delta)$ is:
\begin{align*}
\eps^2 \Delta_g u &= -\eps H_t \dot{\bg}_{\eps} + \ddot{\bg}_{\eps} + \eps^2[\Delta_t \phi - H_t \phi_t + \phi_{tt}] \\
W'(u) &= W'(\bg_{\eps}) + W''(\bg_{\eps}) \phi + Q_0(\phi^2)  \\
\implies 0 &= \eps^2 \Delta_g u - W'(u) \\
&= -\eps H_t \dot{\bg}_{\eps} + L_{\eps}(\phi) + Q_0(\phi^2) + \Rwe \\
&= - \eps H_t \dot{\bg}_{\eps} + L_{\eps}(\phi) + \tilde{Q}_0
\end{align*}
For
\begin{align*}
Q_0(\phi) &= -\frac{1}{2}W'''(\bg_{\eps}) \phi^2 - \phi^3 = \tilde{O}(\phi^2) \\
\tilde{Q}_0(\phi) &= Q_0(\phi) + \Rwe = \tilde{O}(\phi^2, \eps^{\omega}, \eps^{\omega} \phi) \\
&= \tilde{O}(\phi^2, \eps^{\omega})
\end{align*}
\qed
\subsection{Invertibility of the linearized operator, $L_{\eps}$}
In this section, we now prove \cref{SchauderBoundLemma} for the differential operator, $L_{\eps}$: 
\SchauderBoundLemma*
%
\begin{comment}
\begin{lemma} \label{ModifiedPacardBound}
Let $Y$ a $C^{2,\alpha}$ surface and suppose $f: M \to \R$ in $C^{2,\alpha}_{\eps}(M)$ satisfies $f(s,0) \equiv 0$. There exists $\eps_0 > 0$ such that $\forall \eps < \eps_0$, 
\[
||f||_{C^{2,\alpha}_{\eps}(M)} \leq K ||L_{\eps} f||_{C^{\alpha}_{\eps}(M)}
\]
for $K$ independent of $\eps$.
\end{lemma}
\end{comment}
%%
%
\noindent \Pf Rescale the metric to $g_{\eps}$ as in \eqref{rescaledMetric} so that 
\[
L_{\eps} = \eps^2 \Delta_g - W''(\bg_{\eps}) = \Delta_{g_{\eps}} - W''(\bg_{\eps})
\] 
For any $U$ in the interior of $(M, g_{\eps})$, i.e. $\text{dist}(U, Y) > \delta$ fixed, we have that 
\[
||f||_{C^{2,\alpha}_{\eps}(U)} \leq K [||L_{\eps} \phi||_{C^{\alpha}_{\eps}(U)} + ||\phi||_{C^0(U)}]
\]
by schauder theory, with $K$ independent of $\eps$. For points in $Y \times [0, \delta)$, we consider scaled fermi coordinates $\sigma = \eps^{-1} s$, $\tau = \eps^{-1} t$ (along with \eqref{BlowUpMaps}). This gives
\[
||f(\sigma, \tau)||_{C^{2,\alpha}(Y \times [0, \delta \eps^{-1}))} \leq K (||L_{\eps} f(\sigma, \tau)||_{C^{\alpha}(Y \times [0, \delta \eps^{-1}))} + ||f||_{C^0(Y \times [0, \delta \eps^{-1}))})
\]
Note that we've changed $C^{k,\alpha}_{\eps} \to C^{k,\alpha}$ by parameterizing by $(\sigma, \tau)$ instead of $(s,t)$. Moreover, $K$ is independent of $\eps$ by the expansion of $\eps^2 \Delta_g = \Delta_{g_{\eps}}$ in the scaled coordinates, $(\sigma, \tau)$. Undoing the scaling and using compactness of $Y$ and $M^+$, our two bounds give 
\[
||f||_{C^{2,\alpha}_{\eps}(M)} \leq K [||L_{\eps} f||_{C^{\alpha}_{\eps}(M)} + ||f||_{C^0(M)}]
\]
It suffices to prove
\[
||f||_{C^0(M)} \leq \tilde{K} ||L_{\eps} f||_{C^{\alpha}_{\eps}(M)}
\]
for some $\tilde{K} > 0$ also independent of $\eps$. Suppose not, then there exists a sequence of $\{f_j\}$ and $\{\eps_j\}$ such that 
\[
||f_j||_{C^0(M)} \geq j ||L_{\eps} f_j||_{C^{\alpha}_{\eps}(M)}
\]
normalize each $f_j$ by $||f_j||_{C^0}$ so that 
\[
j^{-1} \geq ||L_{\eps} f_j||_{C^{\alpha}_{\eps}(M)} 
\]
Choose $p_j \in M$ so that $|f_j(p_j)| = 1$. By doing a maximum principle comparison with $f \equiv 1$, we see that $\text{dist}(Y,p_j) < \kappa \eps_j$ for some $\kappa > 0$ independent of $\eps_j$ (see appendix \S \ref{TauBounded}). Thus $p_j = (s_j, t_j)$ with $t_j < \kappa \eps_j$. Define the blow ups of $f_j$ around $s_j$:
\begin{align*}
\tilde{f}_j &: B_{\eps^{-1}}(0) \times [0, \delta \eps^{-1}) \to \R  \\
\tilde{f}_j(\sigma, \tau) &= f_j( \eps \sigma + s_j, \eps \tau)
\end{align*}
In $(\sigma, \tau)$ coordinates, we have the local estimate $||L_{\eps} \tilde{f}_j||_{C^{\alpha}(B_{\eps^{-1}}(0) \times [0, \delta \eps^{-1})} \to 0$ since 
\[
||L_{\eps} f_j||_{C^{\alpha}_{\eps}(M)} > ||L_{\eps} \tilde{f}_j||_{C^{\alpha}(B_{\eps^{-1}}(0) \times [0, \delta \eps^{-1})}
\]
Having normalized by $||f_j||_{C^0}$, we have uniform $C^{2,\alpha}_{\eps_j}$ bounds. Moreover $g_{\eps} \to g_{\R^{n}}$ uniformly in $\eps$ after pulling back to $(\sigma, \tau)$ coordinates. Thus, we have uniform $C^{2,\alpha}$ estimates on $\tilde{f}_j(\sigma, \tau)$. Use Arzel\`a Ascoli to pass to a subsequence which converges in $C^2$ to $f_{\infty}: \R^n \times \R^+$ on compact sets. The subsequence comes with $\{\tau_j\} \to \tau^*$ with $0 < \tau^* < \infty$ so that 
\[
|f_{\infty}(0, \tau^*)| = 1
\]
We also have convergence at the boundary, i.e. 
\[
\forall \sigma \in \R^n, \qquad f_{\infty}(\sigma, 0) \equiv 0
\]
This follows because for $C > \tau > 0$, all of the $\tilde{f}_j$'s satisfy 
\[
|\tilde{f}_j(\sigma,\tau)| \leq 2K \tau
\]
by the $C^{2,\alpha}$ bounds and that $\tilde{f}_j(\sigma,0) \equiv 0$. Thus we get the same interior bound for $f_{\infty}$, which forces the same Dirichlet data. Since $j \to \infty \implies \eps_j \to 0$, we get 
\begin{align*}
L_{\eps_j}&= \Delta_{g_{\eps}} - W''(\bg(\tau)) \xrightarrow{j \to \infty} \Delta_{\R^n} + \partial_{\tau}^2 - W''(g) =:\tilde{L} \\
\implies \tilde{L} f_{\infty} &= 0
\end{align*}
this tells us that $f_{\infty}(\sigma, \tau) \equiv 0$ by lemma \ref{HalfSpace} for $\tilde{L} = \Delta_{\R^n} + \partial_{\tau}^2 - W''(g)$ on the half space. But we've point picked so that $f_{\infty}(0, \tau^*) \neq 0$, a contradiction. \qed \nl \nl
From the lemma, we immediately have
\begin{corollary}
For all $\eps < \eps_0$ sufficiently small, the operator
\[
L_{\eps} : C^{2,\alpha}_{\eps, 0}(M) \to C^{\alpha}_{\eps}(M)
\]
is invertible.
\end{corollary}
\noindent In the exact same way, we prove the corresponding $C^{1,\alpha}_{\eps}$ estimate.
%
%Don't think I use this anywhere but whatever it's good to note
%
\begin{lemma} \label{C1ModifiedPacardBound}
Let $Y$ be $C^{2,\alpha}$ and suppose $f: M \to \R$ in $C^{2,\alpha}_{\eps}(M)$ and satisfies $f(s,0) \equiv 0$. There exists $\eps_0 > 0$ such that $\forall \eps < \eps_0$, 
\[
||f||_{C^{1,\alpha}_{\eps}(M)} \leq K ||L_{\eps} f||_{C^0(M)}
\]
\end{lemma}
\noindent \rmk \; Of course, our results apply to $u_{\eps, Y}^{\pm}$ in the closed setting, recreating the lemmas on each of $M^{\pm}$ in the decomposition of $M = M^+ \sqcup_Y M^-$. We contrast the above bounds with the analogous bound in Pacard--Ritore [\cite{pacard2003constant}, Prop 8.6]. In this paper, $\phi: Y \times \R \to \R$ and the bound
\[
||\phi||_{C^{2,\alpha}_{\eps}(Y \times \R)} \leq K ||L_{\eps} \phi||_{C^{\alpha}_{\eps}(Y \times \R)}
\]
holds via the orthogonality condition: $\int_{\R} \phi \dot{g} = 0$, i.e. $\phi$ has been projected away from the kernel of $L_{\eps}$, which allows the authors to exclude $||\cdot ||_{C^0}$ term in the Schauder estimate. Instead of an orthogonality relationship, we use the Dirichlet condition of $u(s,0) \equiv 0 \implies \phi(s,0) \equiv 0$ to get rid of the $||\cdot||_{C^0}$ term.
\begin{comment}
Note that the Dirichlet condition does not imply orthogonality on the half space, as the solution $w(t)$ from \eqref{wEquation} satisfies the Dirichlet condition and also
\[
\int_0^{\infty} w(t) \dot{g}(t) \neq 0
\]
purely by a sign argument.
(\textcolor{purple}{In retrospect, this comment is stupid since we could just take a bump function which vanishes at $t = 0$ or something})
\end{comment}
%Using the exact same proof, we can also get $C^{1,\alpha}_{\eps}$ bounds:
%
%% I don't think we need this lemma, so I'm just commenting it out
%
%
%\begin{corollary} \label{ModifiedPacardBoundC1}
%There $\exists \eps_0 > 0$ such that $\forall \eps < \eps_0$, we have the following: if $\phi$ satisfies $\phi(s,0) \equiv 0$, then
%\[
%||\phi||_{C^{1,\alpha}_{\eps}(B_1(0) \times [-\omega \ln(\eps)))} \leq K ||L_{\eps} \phi||_{C^{-1, \alpha}_{\eps}(B_1(0) \times [0, -\omega\ln(\eps)))}
%\]
%\end{corollary}
%
%
\begin{corollary} \label{InitialEstimateCor}
For $\phi$ satisfying the same conditions as in \cref{SchauderBoundLemma}, we have 
\begin{equation} \label{InitialEstimate}
||\phi||_{C^{2,\alpha}_{\eps}(M)} = O(\eps)
\end{equation}
\end{corollary}
\noindent \Pf \cref{SchauderBoundLemma} and \ref{ACDecomposition} let us conclude that for $u_{\eps}(s,t) = \bg_{\eps}(t) + \phi(s,t)$,
\begin{align*}
||\phi||_{C^{2,\alpha}_{\eps}(M)} &\leq || \eps H_t \dot{\bg}_{\eps}||_{C^{\alpha}_{\eps}(M)} + ||\tilde{Q}_0(\phi)||_{C^{\alpha}_{\eps}(M)} \\
||\tilde{Q}_0(\phi)||_{C^{\alpha}_{\eps}(M)} & \leq ||\phi||_{C^{\alpha}_{\eps}(M)}^2 + K \eps^{\omega} \\
\implies ||\phi||_{C^{2,\alpha}_{\eps}(M)} &\leq O(\eps) + \mu ||\phi||_{C^{\alpha}_{\eps}(M)} \\
||\phi||_{C^{2,\alpha}_{\eps}(M)} &\leq O(\eps)
\end{align*}
where $\mu$ can be made arbitrarily small as $\eps \to 0$ by lemma \ref{InitialPhiBound}.
\subsection{Better tangential behavior}
\label{Tangential}
In this section, we get improved horizontal estimates when $Y$ is $C^{3,\alpha}$. Let $\nabla^Y$ denote the gradient on $Y$, extended as an operator on functions on $Y \times [0, \delta_0)$ in $(s,t)$ coordinates via 
\[
\nabla^Y f(s,t) = g^{ij}(s,0) \partial_{s_i}(f) \partial_{s_j}\Big|_{(s,t)}
\]
\noindent where $\partial_{s_i}$ is identified with $F_*(\partial_{s_i}) \Big|_{s, t}$ using equation \eqref{FMapGeneral}. Then we have
\begin{lemma}\label{ImprovedHorizontal}
Let $Y$ be $C^{3,\alpha}$. For $\phi(s,t)$ as in $u(s,t) = g_{\eps}(t) + \phi(s,t)$ from \ref{InitialDecomposition} and any $\delta > 0$, there exists $K = K(\delta)$ such that
\begin{align} \label{BetterTangent}
||\nabla_Y \phi||_{C^{2,\alpha}_{\eps}(Y \times [0, \delta))} & \leq K \eps
%||D^2_Y \phi||_{C^{2,\alpha}_{\eps}} &= O(\eps)
\end{align}
%
%(\textcolor{purple}{NOte this lemma is local which might screw up changing the notation to all be global on $M$. However we know that $\phi$ has small derivatives inside of $M$, so it all should be fine, just a pain in the ass to write up }) (\textcolor{blue}{- Yes see the lemma where $\phi$ decays exponentially away from $Y$. Though $\nabla^Y$ is no longer defined far away from $Y$ (i.e. fermi coordinates break down unless we're cylindrical), the gradient of $\phi$ is still small!})
\end{lemma}
\noindent \rmk \; The reader may ask how this estimate is ``improved" since $\phi$ satisfies the same $C^{2,\alpha}_{\eps}$ bound. The point is that because of the weighting in the $||\phi||_{C^{2,\alpha}_{\eps}}$ bound, a priori we have
\begin{align*}
||\phi_{s_i}||_{C^{1,\alpha}_{\eps}(Y \times [0, \delta))} &\leq \eps^{-1} ||\phi||_{C^{2,\alpha}_{\eps}(M)} = O(1) \\
||\phi_{t}||_{C^{1,\alpha}_{\eps}(Y \times [0, \delta))} &\leq \eps^{-1} ||\phi||_{C^{2,\alpha}_{\eps}(M)} = O(1)
\end{align*}
by definition of the $C^{k,\alpha}_{\eps}$ norms and \eqref{InitialEstimate}. By contrast, the above lemma \ref{ImprovedHorizontal} gives an $O(\eps)$ bound for $\phi_{s_i}$ near the boundary, i.e. one order in $\epsilon$ better. This method does not work for $\phi_t$ since $[\partial_t, L_{\eps}](\phi)$ is large a priori. However, we do note that for $t > -\omega \eps \ln(\eps)$, the same proof in lemma \ref{InitialPhiBound} gives 
\begin{equation} \label{decayAway}
||\phi||_{C^{k,\alpha}_{\eps}(t > -\omega \eps \ln(\eps))} \leq C(k) \eps^{\omega}
\end{equation}
for $\eps < \eps_0(k)$ and $C(k)$ independent of $\eps$. This will be used below \nl \nl
\noindent \Pf Starting with:
\[
\eps H_t \dot{\bg}_{\eps} (t) = L_{\eps}(\phi) + Q_0(\phi) + \Rwe
\]
move into Fermi coordinates, and apply $\chi_{\delta}(t)\partial_{s_i}$ to each side:
\begin{align*}
\eps (\chi_{\delta}\partial_{s_i} H_t) \dot{\bg}_{\eps} &= L_{\eps}(\chi_{\delta} \phi_{s_i}) + T_0(\phi) \chi_{\delta}\phi_{s_i} + [L_{\eps}, \chi_{\delta} \partial_{s_i}](\phi) + \chi_{\delta}\Rwe \\
&= L_{\eps}(\chi_{\delta}\phi_{s_i}) + E
\end{align*}
where 
\begin{align*}
E &= T_0(\phi) \chi_{\delta}\phi_{s_i} + [L_{\eps}, \chi_{\delta}\partial_{s_i}](\phi) + \chi_{\delta}\Rwe \\
T_0(\phi) &= -\frac{1}{2} W'''(\bg_{\eps}) \phi - 3 \phi^2
\end{align*}
and we can bound
\begin{align*}
||\eps (\chi_{\delta} \partial_{s_i}(H_t) \dot{\bg}_{\eps})||_{C^{\alpha}_{\eps}(M)} & \leq K \eps	
\end{align*}
for the error term, we have
\begin{align*}
||E||_{C^{\alpha}_{\eps}(M)} &\leq ||\chi_{\delta}\phi_{s_i}||_{C^{\alpha}_{\eps}(M)} \cdot ||T_0(\phi)||_{C^{\alpha}_{\eps}(M)} + ||[L_{\eps}, \chi_{\delta}\partial_{s_i}](\phi)||_{C^{\alpha}_{\eps}(M)} + ||\chi_{\delta} \Rwe||_{C^{\alpha}_{\eps}(M)} \\
T_0(\phi) &= \tilde{O}(\phi) \leq K \eps \\
||\chi_{\delta} \Rwe||_{C^{\alpha}_{\eps}(M)}& = o(\eps) \\
\implies ||E||_{C^{\alpha}_{\eps}(M)} & \leq K \eps ||\phi_{s_i}||_{C^{\alpha}_{\eps}(M)} + ||[L_{\eps}, \chi_{\delta}\partial_{s_i}](\phi)||_{C^{\alpha}_{\eps}(M)} + o(\eps) \\
& \leq K \eps + ||[L_{\eps}, \chi_{\delta}\partial_{s_i}](\phi)||_{C^{\alpha}_{\eps}(M)} 
\end{align*}
%
%O(||\phi_{s_i}||_{C^{\alpha}_{\eps}} \cdot ||\phi||_{C^{\alpha}_{\eps}}, \eps) = O(\eps ||\phi_{s_i}||_{C^{\alpha}_{\eps}}) + O(\eps)
%
Here, we've noted that $T_0(\phi) = \tilde{O}(\phi)$ and that $||\phi_{s_i}||_{C^{\alpha}_{\eps}(M)} = O(1)$ a priori. \nl \nl
We now compute the commutator 
\begin{align*}
[L_{\eps}, \chi_{\delta} \partial_{s_i}] &= \eps^2[\Delta_t - H_t \partial_t + \partial_t^2 - W''(\bg_{\eps}), \chi_{\delta}(t) \partial_{s_i}] \\
&= [\eps^2\Delta_t, \chi_{\delta}(t) \partial_{s_i}] + [-\eps^2 H_t \partial_t, \chi_{\delta}(t) \partial_{s_i}] + [\eps^2\partial_t^2, \chi_{\delta}(t) \partial_{s_i}] - [\eps^2 W''(\bg_{\eps}), \chi_{\delta}(t) \partial_{s_i}] 
\end{align*}
in pieces, we have
\begin{align*}
||[\eps^2\Delta_t, \chi_{\delta}(t)\partial_{s_i}] \phi||_{C^{\alpha}_{\eps}(M)} &\leq \eps^2 ||[\Delta_t, \chi_{\delta}(t) \nabla^Y] \phi||_{C^{\alpha}_{\eps}(M)} \\
&\leq \eps^2 ||[\Delta_Y, \nabla^Y] \chi_{\delta}(t)\phi||_{C^{\alpha}_{\eps}(M)} + \eps^2 ||[(\Delta_t - \Delta_Y), \nabla^Y] \chi_{\delta}(t) \phi||_{C^{\alpha}_{\eps}(M)} \\
& \leq \eps^2 ||\chi_{\delta}(t) \Ric(\nabla^Y \phi, \cdot)||_{C^{\alpha}_{\eps}(M)} + \eps^2 ||[(\Delta_t - \Delta_Y), \nabla^Y] \chi_{\delta}(t) \phi||_{C^{\alpha}_{\eps}(M)} \\
& \leq \eps^2 ||\chi_{\delta}(t) \Ric(\nabla^Y \phi, \cdot)||_{C^{\alpha}_{\eps}(M)} + \eps^2 ||[(\Delta_t - \Delta_Y), \nabla^Y] \chi_{\delta}(t) \phi||_{C^{\alpha}_{\eps}(Y \times [0, 2\delta))} \\
%& \leq \eps^2 ||\chi_{\delta}(t) \Ric(\nabla^Y \phi, \cdot)||_{C^{\alpha}_{\eps}(M)} + \eps^2 ||[(\Delta_t - \Delta_Y), \nabla^Y] \chi_{\delta}(t) \phi||_{C^{\alpha}_{\eps}(Y \times [0, -\omega \eps \ln(\eps)))} \\
%& + \eps^2 ||[(\Delta_t - \Delta_Y), \nabla^Y] \chi_{\delta}(t) \phi||_{C^{\alpha}_{\eps}(t > - 2\omega \eps \ln(\eps))} \\
& \leq K \eps  ||\phi||_{C^{1,\alpha}_{\eps}(M)} + K \delta ||\phi||_{C^{2,\alpha}_{\eps}(M)} + O(\eps^{\omega}) \\
& \leq K \delta \eps
\end{align*}
having used \eqref{InitialEstimateCor} and for $K$ independent of $\delta$. We also compute
\begin{align*}
||[-\eps^2 H_t \partial_t, \chi_{\delta}(t) \partial_{s_i}] \phi||_{C^{\alpha}_{\eps}(M)} & \leq \eps^2 [||\chi_{\delta} \partial_{s_i}(H_t) \phi_t||_{C^{\alpha}_{\eps}(M)} + ||H_t \chi_{\delta}'(t) \phi_s||_{C^{\alpha}_{\eps}(M)}] \\
& \leq K \eps^2 ||\phi_t||_{C^{\alpha}_{\eps}(M)} + K \eps^2 \delta^{-1} ||\phi_s||_{C^{\alpha}_{\eps}(M)} \\
& \leq K \delta^{-1} \eps ||\phi||_{C^{1,\alpha}_{\eps}(M)} \\
& \leq K \delta^{-1} \eps^2 
\end{align*}
Furthermore
\begin{align*}
||[\eps^2 \partial_t^2, \chi_{\delta} \partial_{s_i}] \phi||_{C^{\alpha}_{\eps}(M)} &= \eps^2||\chi_{\delta}' \phi_{st} + \chi_{\delta}'' \phi_s||_{C^{\alpha}_{\eps}(M)} \\
& \leq \delta^{-2} ||\phi||_{C^{2,\alpha}_{\eps}(t > \delta)} \\
& \leq K \eps^{\omega} \\
||[\chi_{\delta} \partial_{s_i}, \eps^2 W''(\bg_{\eps})] \phi|| & = 0
\end{align*}
So in conclusion, we have 
\[
||E||_{C^{\alpha}_{\eps}(M)} \leq K \eps
\]
And so 
\begin{align*}
L_{\eps}(\chi_{\delta} \phi_{s_i}) &= \eps (\chi_{\delta} \partial_{s_i}(H_t) \dot{\bg}_{\eps}) - E \\
||\phi_{s_i}||_{C^{2,\alpha}_{\eps}(Y \times [0, \delta))} & \leq ||\chi_{\delta}\phi_{s_i}||_{C^{2,\alpha}_{\eps}(M)} \\
& \leq K ||L_{\eps}(\chi_{\delta} \phi_{s_i})||_{C^{\alpha}_{\eps}(M)} \\
& \leq K \eps
\end{align*}
here, we've used \cref{SchauderBoundLemma} (note that $ \phi_{s_i}(s,0) \equiv 0$). The bound on $\nabla^Y \phi$ now follows. \qed 
\subsection{Proof of theorem \ref{NeumannResult}} \label{ProofOfNeumannResult}
Referencing the decomposition in equation \eqref{ACDecomposition}, we multiply by $\dot{\bg}_{\eps}$ and integrate from $t = 0 $ to  $t = - \omega \eps \ln(\eps)$, picking up a boundary term:
\begin{align} \nonumber 
\tilde{O}(\phi^2) &= -\eps H_t \dot{\bg}_{\eps} + L_{\eps}(\phi) \\ \label{PhiIntegration}
\implies \int_0^{- \omega \eps \ln(\eps)} \tilde{O}(\phi^2) \dot{\bg}_{\eps} &= -\eps \int_0^{- \omega \eps \ln(\eps)} H_t \dot{\bg}_{\eps}^2 + \eps^2 \int_0^{-\omega \eps \ln(\eps)}(\Delta_t \phi) \dot{\bg}_{\eps}  \\ \nonumber
&- \eps^2 \int_0^{- \omega \eps \ln(\eps)} H_t \phi_t \dot{\bg}_{\eps} + \eps^2 \int_0^{- \omega \eps \ln(\eps)} [\phi_{tt} - W''(\bg_{\eps}) \phi] \dot{\bg}_{\eps}
\end{align}
Note that the left hand side of \eqref{PhiIntegration} can be bounded
\begin{align*}
\Big|\Big|\int_0^{- \omega \eps \ln(\eps)} \tilde{O}(\phi^2) \dot{\bg}_{\eps} dt\Big|\Big|_{C^{\alpha}(Y)} & \leq K \eps^{2-\alpha} \int_0^{- \omega\eps \ln(\eps)} \dot{\bg}_{\eps} dt \\
& \leq K \eps^{3-\alpha} + O(\eps^\omega) \\
& \leq K \eps^{3-\alpha}
\end{align*}
since $\omega > 5$. For the right hand side of \eqref{PhiIntegration}, we note
\[
|H_t - H_0| \leq t ||\dot{H}_t||_{C^0(Y \times [0, \delta))} \leq K t
\]
which follows by the mean value theorem and \eqref{secondFundamentalExpansion}. Moreover
\[
\implies \int_0^{- \omega\eps \ln(\eps)} \dot{\bg}_{\eps}^2 = \eps \sigma_0 + C \eps^{\omega}
\]
We further note from \eqref{BetterTangent}:
\begin{align*}
\Big| \Big| \int_0^{-\omega\eps \ln(\eps) } \eps^2 (\Delta_t \phi) \dot{\bg}_{\eps} \Big| \Big|_{C^{\alpha}(Y)} & = \Big| \Big| \int_0^{- \omega\eps \ln(\eps) } \eps^2 (\Delta_0 \phi) \dot{\bg}_{\eps} \Big| \Big|_{C^{\alpha}_{Y}} + \Big| \Big| \int_0^{- \omega\eps \ln(\eps) } \eps^2 ([\Delta_t - \Delta_0]\phi) \dot{\bg}_{\eps} \Big| \Big|_{C^{\alpha}_{Y}} \\
& = O(\eps^{3 - \alpha}) + O(\eps^{4-\alpha}) \\
&= O(\eps^{3-\alpha})
\end{align*}
Here, we've again used that $\Delta_t - \Delta_0 = t L$ where $L$ is a second order linear differential operator with bounded coefficients. In both cases, we use equation \eqref{ImprovedHorizontal} to  get the final bounds. Similarly, from equation \eqref{InitialEstimate}, we have
\begin{align*}
\Big| \Big| \int_0^{- \omega \eps \ln(\eps)} \eps^2 H_t \phi_t \dot{\bg}_{\eps} \Big| \Big|_{C^{\alpha}(Y)} &= O(\eps^{3 - \alpha}) 
\end{align*}
We also compute
\begin{align*}
\int_{0}^{- \omega\eps \ln(\eps)} \dot{\bg}_{\eps} &= \eps(1 + O(\eps^\omega)) \\
\int_0^{- \omega\eps \ln(\eps)} \left( \eps^2 \phi_{tt} \dot{\bg}_{\eps} - W''(\bg_{\eps}) \phi \dot{\bg}_{\eps} \right) &= -\eps^2 \sigma \phi_t(s,0) + O(\eps^k) \\
\int_0^{-\omega\eps \ln(\eps)} \eps H_t \dot{\bg}_{\eps}^2	&= \eps^2 H_0 \sigma_0 + O(\eps^3)
\end{align*}
with all estimates holding in $C^{\alpha}(Y)$. Combining and noting that $\sigma_0 \sigma^{-1} = \frac{2}{3}$, we have
\[
\eps^2 \phi_t(s,0) = - \eps^2 \frac{2}{3} H_0  + O(\eps^{3-\alpha})
\]
for any $k \geq 3$. We summarize this as
\[
\phi_t(s,0) = -\frac{2}{3} H_0  + O(\eps^{1-\alpha})
\]
so that 
\begin{equation} \label{NormalDerivative}
\partial_{\nu} u(p) = \partial_t u = \eps^{-1}\dot{g}(0) + \phi_t(p,0) = \frac{1}{\eps \sqrt{2}}  - \frac{2}{3} H_Y(p) + O(\eps^{1-\alpha}) 
\end{equation}
holds in $C^{\alpha}(Y)$. This proves \cref{NeumannResult}. \qed

\section{Higher Order Expansions} \label{HigherOrderSection}

\subsection{Next Order Expansion for $H_Y = 0$} \label{NextOrderExpansionMinimal}
In this section, we give a more precise description of the normal derivative for $Y$ minimal and prove \cref{ImprovedNeumannMinimal}:
\ImprovedNeumannMinimal*
\noindent \Pf When $Y$ is minimal, (and hence smooth since we've assumed all our hypersurfaces are at least $C^{2,\alpha}$), \eqref{ACDecomposition} becomes
\[
L_{\eps}(\phi) = \eps H_t \dot{\bg}_{\eps} + Q_0(\phi^2)
\]
if $H_Y = H_0 = 0$, then 
\[
H_t = \int_0^t \dot{H}_r dr
\]
we expand this as
\begin{equation} \label{PhiExpansionForMinimal}
L_{\eps}(\phi) = \eps^2 \dot{H}_0(s) \left( \frac{t}{\eps} \right)\dot{\bg}_{\eps}(t) + \eps \left(\int_0^t \int_0^r \ddot{H}_w(s) dw dr\right) \dot{\bg}_{\eps}(t) + Q_0(\phi^2)
\end{equation}
with the goal of showing 
\[
\Big| \Big| \eps \left(\int_0^t \int_0^r \ddot{H}_w(s) dw dr\right) \dot{\bg}_{\eps}(t) \Big| \Big|_{C^{\alpha}_{\eps}(M)} = o(\eps^2)
\]
The $C^0$ bound holds clearly as 
\begin{align*}
\Big| \eps \left(\int_0^t \int_0^r \ddot{H}_w dw dr\right) \dot{\bg}_{\eps}(t) \Big| & \leq \eps \left(\int_0^t \int_0^r \sup_{\substack{ s \in Y \\ w \in [0, - \omega \eps \ln(\eps))}}|\ddot{H}_w(s)| dw dr \right) \dot{\bg}_{\eps} \\
& \leq K \eps t^2 \dot{\bg}_{\eps} \\
&\leq K \eps^3  \sup_{t \in [0, - \omega \eps \ln(\eps))}\Big|\left(\frac{t}{\eps} \right)^2 \dot{\bg}_{\eps}(t) \Big| \\
& \leq K \eps^3
\end{align*}
For the $[\cdot]_{\alpha}$ bound we have 
\begin{align*}
f(s,t) &:= \eps\left(\int_0^t \int_0^r \ddot{H}_w(s) dw dr\right) \dot{\bg}_{\eps}(t) \\
[f]_{\alpha, Y \times [0, - 2\omega \eps \ln(\eps))} &\leq \eps \Big( \Big| \Big|\int_0^t \int_0^r \ddot{H}_w(s) dw dr \Big| \Big|_{C^0(Y \times [0, - 2\omega \eps \ln(\eps))} [\dot{\bg}_{\eps}]_{\alpha, Y \times [0, - 2\omega \eps \ln(\eps)}  \\
& \qquad \quad + \Big[\int_0^t \int_0^r \ddot{H}_w(s) dw dr \Big]_{\alpha, Y \times [0, - 2\omega \eps \ln(\eps)} ||\dot{\bg}_{\eps}||_{C^0(Y \times [0, - 2\omega \eps \ln(\eps))}\Big)
\end{align*}
having noted that $\dot{\bg}_{\eps} \equiv 0$ for $t > -2 \omega \eps \ln(\eps)$. On $Y \times [0, - 2\omega \eps \ln(\eps))$, these norms are bounded by
\begin{align*}
\Big| \Big|\int_0^t \int_0^r \ddot{H}_w(s) dw dr \Big| \Big|_{C^0(Y \times [0, - 2\omega \eps \ln(\eps))} & \leq K t^2 = O(\eps^2 \ln(\eps)^2) \\
[\int_0^t \int_0^r \ddot{H}_w(s) dw dr ]_{\alpha, Y \times [0, - 2\omega \eps \ln(\eps)} & \leq K t^{2-\alpha} = O(\eps^{2-\alpha} \ln(\eps)^{2-\alpha}) \\
||\dot{\bg}_{\eps}||_{C^0(Y \times [0, - 2\omega \eps \ln(\eps))} &= O(1) \\
[\dot{\bg}_{\eps}]_{\alpha, Y \times [0, - 2\omega \eps \ln(\eps)} &= O(\eps^{-\alpha})
\end{align*}
so that 
\begin{align*}
[f]_{\alpha, \eps} & \leq O(\eps^{3-2\alpha} \ln(\eps)^{2-\alpha}) \\
&= o(\eps^2) 
\end{align*}
Further noting that
\[
||Q_0(\phi)||_{C^{\alpha}_{\eps}(M)} \leq K ||\phi||_{C^{\alpha}(M)} + K\eps^{\omega} \leq K \eps^2
\]
We then have to leading order
\[
L_{\eps}(\phi) = O(\eps^2)
\]
in $C^{\alpha}_{\eps}(M)$. From \cref{SchauderBoundLemma},
\[
||\phi||_{C^{2,\alpha}_{\eps}(M)} \leq K \eps^{2}
\]
If we differentiate \eqref{PhiExpansionForMinimal} with respect to $\partial_{s_i}$ again, we get by the same bounding techniques
\begin{equation} \label{TangentialDerivInterEq}
L_{\eps}(\partial_{s_i} \phi) = \eps \left(\int_0^{t} \partial_{s_i}(\dot{H}_r) dr \right) \dot{\bg}_{\eps} + \overline{Q}_0(\phi \phi_{s_i}, \eps^2 D^2\phi, \eps^2 D \phi) + o(\eps^2)
\end{equation}
Using our bound on $||\phi||_{C^{2,\alpha}_{\eps}}$ and \cref{SchauderBoundLemma} composed with $\chi_{\delta}$ as in lemma \ref{ImprovedHorizontal}, we get
\[
||\phi_{s_i}||_{C^{2,\alpha}_{\eps}(Y \times [0, \delta))} = O(\eps^2)
\]
so that
\[
\sup_{t \in [0, - \omega \eps \ln(\eps))} ||\Delta_t \phi(s,t)||_{C^{\alpha}(Y)} \leq ||\Delta_t \phi||_{C^{\alpha}(Y \times [0, \delta))} = O(\eps^{1 - \alpha})
\]
%Remember that this is the right tangential estimate (it should only increase if we did $\partial_{\sigma_i} = \eps \partial_{s_i}$
%
Now we multiply \eqref{PhiExpansionForMinimal} by $\dot{\bg}_{\eps}$ and integrate
\begin{align*}
\int_0^{-\omega \eps \ln(\eps)} L_{\eps}(\phi) \dot{\bg}_{\eps}(t) dt &= \int_0^{- \omega\eps \ln(\eps)} \left( \eps^2[\Delta_t(\phi) - H_t \phi_t + \phi_{tt}] - W''(\bg_{\eps})\phi \right)\dot{\bg}_{\eps}(t) dt \\
&= - \eps^2 \sigma \phi_t(s,0) + O(\eps^{4-\alpha}) \\
\eps^2 \int_0^{- \omega\eps \ln(\eps)} \dot{H}_0(s) \left( \frac{t}{\eps} \right)\dot{\bg}_{\eps}^2(t) &= \kappa_0 \eps^3 \dot{H}_0 \\
\Big|\eps \int_0^{- \omega\eps \ln(\eps)}\left(\int_0^t \int_0^r \ddot{H}_w(s) dw dr\right) \dot{\bg}_{\eps}(t)^2 dt\Big| & \leq C \eps \int_0^{- \omega \eps \ln(\eps)} t^2 \dot{\bg}_{\eps}(t)^2 dt \\
& = O(\eps^4) \\
\int_0^{- \omega\eps \ln(\eps)} Q_0(\phi) \dot{\bg}_{\eps} dt &= O(\eps^4) \\
\end{align*}
where these error terms hold in $C^{\alpha}(Y)$. Now note that
\[
\dot{H}_0 = [\Ric(\nu, \nu) + |A_Y|^2] 
\]
so that 
\begin{equation} \label{NormalDerivativeMinimal}
\phi_t(s,0) = \sigma^{-1}\kappa_0 \eps [\Ric(\nu, \nu) + |A_Y|^2] + O(\eps^{2-\alpha})
\end{equation} \qed 
%
%Note that $\Ric(\nu, \nu) = \Ric_M(\nu, \nu)$ and $|A_Y|^2$ are both \textit{independent of the choice of normal}, i.e. independent of the orientation, so the Neumann difference operator here will actually be $O(\eps^2)$.
%%
\subsection{Full characterization of Neumann Data}
One can compare theorems \ref{NeumannResult} and \ref{ImprovedNeumannMinimal} and note that more terms can be gleaned. In fact, if $Y$ is a $C^{k + 3, \alpha}$ surface, we can find an expansion for $u_{\eps}(s,t)$ (and hence, $\partial_t u_{\eps}|_{t = 0}$) up to order $k$ ($k-1$ respectively). Let
\[
a_{i,j}(s):= a_{i,j}\left(\{\partial_s^{\beta} \partial_t^j H_t |_{t = 0}\}_{j + |\beta| \leq i}\right)(s)
\] 
denote a polynomial in derivatives of $H_t(s)$ at $t = 0$. Define 
\[
i \in \Z^{\geq 0}, \qquad \sigma(i) := \max(0, 2 \lceil i/2 \rceil - 2) = \begin{cases}
0 & i = 0 \\
\text{largest even integer less than } i & i > 0
\end{cases}
\] 
We prove our main result, \cref{InductiveExpansion}:
\InductiveExpansion*
\noindent \Pf \; We actually prove the following by induction: for any $k < \overline{k}$, we have
\begin{align*}
u_{\eps}^+(s,t) &= \bg_{\eps}(t) + \sum_{i = 1}^{k} \eps^i \cdot \left(\sum_{j = 0}^{m_i} a_{i,j}(s) \bw_{i,j}(t/\eps) \right) + \phi \\
L_{\eps}(\phi) &= R_{k+1}(s,t) + F_{k}(\phi) \\
%	L_{\eps}(\phi) &= \sum_{i = 1}^{n_k} \eps^{k + i} \sum_{j = 0}^{N_i} b_{i,j}(s) f_{i,j}(t/\eps) + \tilde{O}(\eps \phi, \phi^2) \\
||a_{i,j}||_{C^{\alpha}(Y)} &= O(1)\\
%	||b_{i,j}||_{C^{\alpha}(Y)} &= O(1) \\
%	||f_{i,j}||_{C^{\alpha}([0, \infty))} &= O(1) \\ 
||w_{i,j}||_{C^{k,\alpha}([0, \infty))} &= O(1)
\end{align*}
Where $\{w_{i,j}(t)\}$ are all exponentially decaying. Moreover, we require that $R_{k+1}$ can be expanded in powers of $\eps$ to arbitrary order less than $\omega$ as follows:
\begin{align*}
\forall \ell \geq 0, \qquad \exists &\{b_{k+1,i,j}(s)\}, \quad \{f_{k+1,i,j}(t)\} \quad \st \\
R_{k+1}(s,t) &= \eps^{k+1} \sum_{j = 0}^{N_{k+1}} b_{k + 1,0,j }(s) \overline{f}_{k + 1, 0, j}(t/\eps)  \\
&+ \sum_{i = 1}^{\ell} \eps^{k + 1 + i} \sum_{j = 0}^{N_{k+1,i}} b_{k+1,i,j}(s) \overline{f}_{k + 1,i,j}(t/\eps) + O(\eps^{\ell + k + 2}) \\
||f_{k+1,i,j}||_{C^{\alpha}([0, \infty))} &= O(1) 
\end{align*}
\noindent where the expansion holds $R_{k+1}$ holds in $C^{\alpha}_{\eps}(M)$ for $\ell + k + 2 \leq \overline{k}+1$ assuming $\overline{k} - k - 1 \geq 0$. In this sense, we see that there is a partial expansion of the remainder up to any order. Here, we require that 
\begin{itemize}
\item $b_{k+1,0,j}(s) = b_{k+1,0,j}\left( \{\partial_t^p \partial_{s}^{\beta} H_t \Big|_{t = 0} \right)$ depends on \textbf{at most $\sigma(k+1)$ tangential derivatives} of $\{\partial_t^p H_t\Big|_{t = 0}(s)\}$. 
\item For $i \geq 1$, $b_{k + 1, i,j}(s) = b_{k + 1, i, j}\left( \{\partial_t^p \partial_{s}^{\beta} H_t \Big|_{t = 0} \right)$ is a polynomial in \textbf{at most $\sigma(k + 2)$ tangential derivatives} of $\partial_t^p H_t\Big|_{t = 0}(s)$. 
\item Each $f_{k+1, i,j}(t)$ is exponentially decaying in $C^{\infty}$ and $\overline{f}_{k+1, i, j}$ is the modification with a smooth cutoff. This allows us to solve
\begin{align} \label{IterativeLinearizedEq} \nonumber
w_{k+1,i,j}&: [0, \infty) \to \R \\
\ddot{w}_{k+1,i,j}(t) - W''(g(t)) w_{k+1,i,j}(t) &= f_{k+1,i,j}(t) \\ \nonumber
w_{k+1,i,j}(0) &= 0 \\
\lim_{t \to \infty} w_{k+1,i,j}(t) & = 0 
\end{align}
by section \S \ref{HalfLineSection} in the appendix. 
\end{itemize}
We also require that $F_{k}(\phi)$ is an error term which has at most cubic dependency on $\phi$ in the following form: 
\begin{align*}
F_k(\phi) &= \eps\left[ \sum_{i = 1}^{n_k} c_{k,i}(s) \overline{h}_{k,i}(t/\eps) \right] \phi \\
& \;\;\; + \left[\sum_{i = 1}^{m_k} d_{k,i}(s) \overline{p}_{k,i}(t/\eps)\right] \phi^2 \\
& \;\;\; - \phi^3 \\
||h_{k,i}(t)||_{C^{\alpha}([0, \infty))} &= O(1) \\
||p_{k,i}(t)||_{C^{\alpha}([0, \infty))} &= O(1)
\end{align*}
Moreover
\begin{itemize}
\item $\{h_{k,i}\}$ and $\{p_{k,i}\}$ are expontentially decaying in $C^{\infty}$
\item $\{c_{k,i}\}$, $\{d_{k,i}\}$ depend on \textbf{at most $\sigma(k)$ tangential derivatives} of $\partial_t^p H_t\Big|_{t = 0}(s)$.
\end{itemize}
Note that $L_{\eps}(\phi) = R_{k+1} + F_{k}$ and \cref{SchauderBoundLemma} automatically gives the conclusion of 
\[
||\phi||_{C^{2,\alpha}_{\eps}(M)} = O(\eps^{k+1})
\]
From hereon in the proof, we assume that $\omega > \overline{k} + 2$. \nl \nl
\noindent \textbf{Base Case $k = 0$}: \nl 
This is the content of corollary \ref{InitialEstimateCor}
\begin{align*}
u_{\eps}^+(s,t) &= \bg_{\eps}(t) + \phi(s,t) \\
L_{\eps}(\phi) &= \left(\eps H_t \dot{\bg}_{\eps} - \Rwe\right) - \left[ \frac{1}{2} W'''(\bg_{\eps}) \phi^2  + \phi^3 \right]  \\
&= R_{1}(s,t) + F_0(\phi)
\end{align*}
from \eqref{ACDecomposition}. At this level of expansion, $\{a_{1,i,j} = 0\}, \; \{w_{1,i,j} = 0\}$. We see that $R_1(s,t) = \eps H_t \dot{\bg}_{\eps} - \Rwe$ satisfies our inductive assumptions simply by expanding
\begin{align*}
H_t &= H_0 + t \dot{H} + \frac{t^2}{2} \ddot{H} + \dots + \frac{t^{\ell + 1}}{(\ell + 1)!} \partial_t^{\ell + 1} H_t |_{t = 0} + O(t^{\ell + 2}) \\
\implies \eps H_t \dot{\bg}_{\eps} &= \eps (H_0) \dot{\bg}_{\eps} + \sum_{i = 1}^{\ell} \eps^{i+1} \left(\frac{1}{(i + 1)!} \partial_t^{\ell + 1} H_t \Big|_{t = 0} \right) \left[\left( \frac{t}{\eps} \right)^{i + 1} \dot{\bg}_{\eps}(t) \right] + O(\eps^{\ell + 2}) \\
\Rwe &= \ddot{\bg}_{\eps} - W'(\bg_{\eps}) = O(\eps^{\omega}) = O(\eps)
\end{align*}
since $\omega > \overline{k} + 2$. Here, we've noted that $\left(\frac{t}{\eps}\right)^i \dot{\bg}_{\eps}$ is bounded in $C^{k,\ell}_{\eps}$ for all $\ell, \; \alpha, \; i$. Thus $R_1(s,t)$ satisfies our inductive assumptions. Note that computing $\partial_t^{i} H_t$ does not require extra regularity of $Y$ - simply expand 
\eqref{secondFundamentalExpansion} in $t$. In particular,
\begin{align*}
N_1 &= N_{1,i} = 0 \\
b_{1,i,0}(s) &= \left(\frac{1}{(i + 1)!} \partial_t^{\ell + 1} H_t \Big|_{t = 0} \right) \\
f_{1,i,0}(s) &= \left( \frac{t}{\eps} \right)^{i + 1} \dot{\bg}_{\eps}(t)
\end{align*}
Moreover, each $b_{i+1}$ depends on $0 \leq \sigma(1), \sigma(2)$ tangential derivatives of $\{\partial_t^{p} H_t |_{t = 0}\} $. And finally, each $f_{1,i, 0}(t)$ is exponentially decaying. Similarly, it is clear that $F_0(\phi)$ satisfies our inductive assumptions as it only has quadratic and cubic terms with bounded coefficients in $C^{k,\alpha}_{\eps}$ that are also exponentially decaying.
\begin{align*}
%c_{0,1}(s) &= 1 \\
%\overline{h}_{1,0}(t/\eps) &= -\frac{1}{\eps} [W''(g_{\eps}) - W''(\bg_{\eps})] \\
d_{0,1}(s) &= 1 \\
p_{0,1}(t/\eps) &= \frac{1}{2} W'''(\bg_{\eps}) = 3 \bg_{\eps}
\end{align*}
\noindent \textbf{Induction}: \nl
Now assume that we have an expansion up to order $k-1$ for $k \leq \overline{k}$:
\begin{align*}
u_{\eps}^+(s,t) &= \bg_{\eps}(t) + \sum_{i = 1}^{k-1} \eps^i \sum_{j = 0}^{m_i} a_{i,j}(s) \bw_{i,j}(t/\eps) + \phi \\
L_{\eps}(\phi) &= R_{k}(s,t) + F_{k-1}(\phi)
\end{align*}
%
%so that $\Rwe = O(\eps^{k+1})$ . Note that $\Rwe$ is exponentially decaying since it is supported on $[-\omega \eps \ln(\eps), -2 \omega \eps \ln(\eps)]$. 
We expand (for any $\ell \geq 1$)
\begin{equation} \label{RkPartialExpansion}
R_k(s,t) = \eps^k \sum_{j = 0}^{N_k} b_{k,0,j}(s) \overline{f}_{k,0,j}(t/\eps) + \left[\sum_{i = 1}^{\ell} \eps^{i + k} \sum_{j = 0}^{N_{k,i}} b_{k,i,j}(s) \overline{f}_{k,i,j}(t/\eps) + O(\eps^{\ell + 1 + k})\right]
\end{equation}
where we know $\{b_{k,0,j}\}$ depend on at most $\sigma(k-1 + 1)= \sigma(k) = \max(0, 2 \lceil k/2 \rceil - 2) \leq k-1$ derivatives of $\partial_t^p H_t\Big|_{t = 0}(s)$. We can compute $2$ more tangential derivatives of $b_{k,0,j}$ when $Y$ is $C^{k + 3, \alpha}$. With this, we use \eqref{RkPartialExpansion} and write
\[
\phi(s,t) = \eps^k \sum_{j = 0}^{N_k} b_{k,0,j}(s) w_{k,0,j}(t/\eps) + \tilde{\phi}
\]
such that each $w_{k,0,j}$ solves
\begin{align*}
\ddot{w}_{k,0,j}(t) - W''(g_{\eps}) w_{k,0,j}(t) &= f_{k,0,j}(t) 
\end{align*}
and has bounded $C^{\alpha}(\R^+)$ norm and is exponentially decaying. This follows again by \S \ref{HalfLineSection}. Multiplying these functions by cutoffs, we get
\[
||\ddot{\overline{w}}_{k,0,j}(t) - W''(\bg_{\eps}) \overline{w}_{k,0,j}(t) - \overline{f}_{k,0,j}(t)||_{C^{\alpha}(\R^+)} \leq C_{k,0,j} \eps^{\omega} \leq C_{k,0,j} \eps^{k+2}
\]
for some constants $C_{k,0,j}$ independent of $\eps$. With this expansion, we have 
\begin{align*}
L_{\eps}(\phi) &= L_{\eps}(\tilde{\phi}) \\
&+ \eps^{k+2} \sum_{j = 0}^{N_{k}} \Delta_t(b_{k,0,j})(s) \bw_{k,0,j}(t/\eps) \\
&- \eps^{k+1} \sum_{j = 0}^{N_{k}} H_t(s) b_{k,0,j}(s) \dot{\bw}_{k,0,j}(t/\eps) \\
& + \eps^k \sum_{j = 0}^{N_{k}} b_{k,0,j}(s) [\ddot{\bw}_{k,0,j}(t/\eps) - W''(\bg_{\eps}) \bw_{k,0,j}(t/\eps)]
\end{align*} 
Because $\{b_{k,0,j}\}$ depend on $\sigma(k) \leq k - 1$ derivatives of $\partial_t^p H_t \Big|_{t = 0}(s)$, we know that $\Delta_t b_{k,0,j}$ is at least in $C^{\alpha}(Y)$ since $k \leq \overline{k}$ and $Y$ is $C^{\overline{k} + 3, \alpha}$. Using \eqref{IterativeLinearizedEq}, we see that the last line cancels with the first term in \eqref{RkPartialExpansion} at the cost of an $O(\eps^{k + 2})$ error. We also expand
\begin{align*}
\Delta_t(b_{k,0,j})\bw_{k,0,j}(t/\eps) &= \Delta_0(b_{k,0,j})(s)\bw_{k,0,j}(t/\eps) + \eps \left( \frac{\Delta_t - \Delta_0}{t} \right)(b_{k,0,j}) \cdot \left(\frac{t}{\eps}\right) \bw_{k,0,j}(t/\eps) \\
H_t(s) b_{k,0,j}(s) \dot{w}_{k,0,j}(t/\eps) &= H_0(s) b_{k,0,j}(s) \dot{\bw}_{k,0,j}(t/\eps) + \eps\left(\frac{H_t(s) - H_0(s)}{t}\right) \left(\frac{t}{\eps}\right) b_{k,0,j}(s) \dot{\bw}_{k,0,j}(t/\eps) 
\end{align*}
where we can write 
\begin{align*}
(\Delta_t b_{k,0,j})&= \sum_{i = 0}^{m} \frac{t^i}{i!} \left(\partial_t^i \Delta_t\Big|_{t = 0}\right) (b_{k,0,j})(s) + O(t^{m + 1}) \\
(\partial_t^m \Delta_t)\Big|_{t = 0} &: C^{l + 2}(Y) \to C^{l} \\
&= \left(\partial_t^m g^{ij}(s, t)\Big|_{t = 0} \right) \partial_{s_i} \partial_{s_j} - \left( \partial_t^m b^k(s,t) \Big|_{t = 0} \right) \partial_{s_k}  
\end{align*}
These expansions in $t$ do not require higher regularity of $H_t(s)$, as can be seen from the expansion of the metric, $g(s,t)$, and the second fundamental form, $A(s,t)$, in equations \eqref{metricExpansion} and \eqref{secondFundamentalExpansion}. This allows us to make sense of $\frac{\Delta_t - \Delta_0}{t}$. Similarly
\begin{align*}
H_t(s) b_{k,0,j}(s) \dot{\bw}_{k,0,j}(t/\eps) &= \sum_{i = 0}^{m} \eps^i \frac{1}{i!} \left( \partial_t^i H_t\Big|_{t = 0}(s) \right) \left[ \left(\frac{t}{\eps}\right)^i \dot{\bw}_{k,0,j}(t/\eps)  \right] + O(\eps^{m + 1})
\end{align*}
for any $m$.
%
\begin{comment}
We can also use the same techniques as in the base case to expand 
%
\begin{align*}
\Delta_t(b_k)w_{k,j}(t/\eps) &= \Delta_0(b_k)(s)w_{k,j}(t/\eps) + \eps \left( \frac{\Delta_t - \Delta_0}{t} \right) \left(\frac{t}{\eps}\right) (b_k) w_{k,j}(t/\eps) \\
H_t(s) b_k(s) \dot{w}_{k,j}(t/\eps) &= H_0(s) b_k(s) \dot{w}_{k,j}(t/\eps) + \eps\left(\frac{H_t(s) - H_0(s)}{t}\right) \left(\frac{t}{\eps}\right) b_k(s) \dot{w}_{k,j}(t/\eps)
\end{align*}
%
\end{comment} 
Similarly, we have
\begin{align*}
F_{k-1}(\phi) &= \eps\left[ \sum_{i = 1}^{n_{k-1}} c_{k-1,i}(s) \bh_{k-1,i}(t/\eps) \right] \phi + \left[\sum_{i = 1}^{m_{k-1}} d_{k-1,i}(s) \overline{p}_{k-1,i}(t/\eps)\right] \phi^2 - \phi^3 \\
&= \eps\left[ \sum_{i = 1}^{n_{k-1}} c_{k-1,i}(s) \bh_{k-1,i}(t/\eps) \right] \left( \eps^k \sum_{j = 0}^{N_k} b_{k,0,j}(s) \overline{f}_{k,0,j}(t/\eps) + \tilde{\phi} \right)\\
& \;\;\; + \left[\sum_{i = 1}^{m_{k-1}} d_{k-1,i}(s) \overline{p}_{k-1,i}(t/\eps)\right] \left( \eps^k \sum_{j = 0}^{N_k} b_{k,0,j}(s) \overline{f}_{k,0,j}(t/\eps) + \tilde{\phi} \right)^2 \\
& \;\;\; - \left( \eps^k \sum_{j = 0}^{N_k} b_{k,0,j}(s) \overline{f}_{k,0,j}(t/\eps) + \tilde{\phi} \right)^3 \\	
\end{align*}
with $\{c_{k-1,i}\}$, $\{d_{k-1,i}\}$ depending on at most $\sigma(k-1)$ derivatives of $\partial_t^p H_t \Big|_{t = 0}(s)$. If we expand and relabel, noting that the product of exponentially decaying functions are themselves exponentially decaying, we get
\begin{align*}
F_{k-1}(\phi) &= \eps^{k+1}\left[ \sum_{i = 1}^{\tilde{n}_k} C_{k,i}(s) \overline{h}^*_{k,i}(t/\eps) \right] \\
& \;\;\; + \eps \left[ \sum_{i = 1}^{n_{k}} c_{k,i}(s) \overline{h}_{k,i}(t/\eps) \right] \tilde{\phi} \\
& \;\;\; + \left[\sum_{i = 1}^{m_{k}} d_{k,i}(s) \overline{p}_{k,i}(t/\eps)\right] \tilde{\phi}^2 \\
& \;\;\; - \tilde{\phi}^3 
\end{align*}
for some $n_k, m_k, \tilde{n}_k$. Here
\begin{itemize}
\item $\{h^*_{k,i}, h_{k,i}, p_{k,i}\}$ are all exponentially decaying and $O(1)$ in $C^{\alpha}$ norm 
\item $\{C_{k,i}\}$, $\{c_{k,i}\}$, $\{d_{k,i}\}$ depend on at most $\sigma(k)$ derivatives of $\partial_t^p H_t \Big|_{t = 0}(s)$. 
\end{itemize} 
We define
\begin{align*}
R_{k+1}(s,t) &:= \left[ R_k(s,t) - \eps^k \sum_{j = 0}^{N_k} b_{k,0,j}(s) \overline{f}_{k,0,j}(t/\eps) \right] + \eps^{k+1}\left[ \sum_{i = 1}^{\tilde{n}_k} C_{k,i}(s) \overline{h}^*_{k,i}(t/\eps) \right]\\
& -  \eps^{k+2} \sum_{j = 0}^{N_{k}} \Delta_t(b_{k,0,j})(s) \overline{w}_{k,0,j}(t/\eps) + \eps^{k+1} \sum_{j = 0}^{N_{k}} H_t(s) b_{k,0,j}(s) \dot{\bw}_{k,0,j}(t/\eps) \\
\implies R_{k+1}(s,t) &= \eps^{k+1} \left[\sum_{j = 0}^{N_{k,1}} b_{k,1,j}(s) \overline{f}_{k,1,j}(t/\eps) + \sum_{i = 1}^{\tilde{n}_k} C_{k,i}(s) \overline{h}^*_{k,i}(t/\eps) + \sum_{j = 0}^{N_{k}} H_0(s) b_{k,0,j}(s) \dot{\bw}_{k,0,j}(t/\eps) \right]  \\
& + \eps^{k+2} \sum_{j = 0}^{N_{k}} \left(\frac{H_t(s) - H_0(s)}{t} \right)  \left(\frac{t}{\eps}\right) b_{k,0,j}(s) \dot{\bw}_{k,0,j}(t/\eps)  -  \eps^{k+2} \sum_{j = 0}^{N_{k}} \Delta_t(b_{k,0,j})(s) \overline{w}_{k,0,j}(t/\eps) \\
& + \sum_{i = 2}^{\ell} \eps^{i + k} \sum_{j = 0}^{N_{k,i}} b_{k,i,j}(s) \overline{f}_{k,i,j}(t/\eps) + O(\eps^{\ell + 1 + k})
\end{align*}
moreover, note that $\{b_{k+1,i,j}\}$ depends on at most $\sigma(k+1)$ derivatives for $i \geq 1$, while $C_{k,i}(s)$ and $b_{k,0,j}(s)$ depend on at most $\sigma(k)$ derivatives. Also recalling that $\ell$ is any value such that $\ell + 1 + k \leq \overline{k} + 1$, we can rewrite the above as 
\begin{align*}
R_{k+1}(s,t) &= \eps^{k+1} \sum_{j = 0}^{N_{k+1}} b_{k+1,0,j}(s) \overline{f}_{k+1,0, j}(t/\eps) \\
& + \sum_{i = 1}^{\ell} \eps^{i + k + 1} \sum_{j = 0}^{N_{k+1,i}} b_{k + 1,i, j}(s) \overline{f}_{k + 1, i, j}(t/\eps) + O(\eps^{\ell + k + 2})
\end{align*}
adjusting the expansion depending on the value of $\overline{k} - k - 1$. If $k = \overline{k}$, then no such expansion is needed, as we've reached the maximal value of $k$ in the induction. Furthermore, we define
\[
F_k(\tilde{\phi}):= \eps \left[ \sum_{i = 1}^{n_k} c_{k,i}(s) \overline{h}_{k,i}(t/\eps) \right] \tilde{\phi} + \left[\sum_{i = 1}^{m_k} d_{k,i}(s) \overline{p}_{k,i}(t/\eps)\right] \tilde{\phi}^2 - \tilde{\phi}^3
\]
So that
\[
L_{\eps}(\tilde{\phi}) = R_{k+1}(s,t) + F_k(\tilde{\phi})
\]
with the correct decomposition and regularity of coefficients. Now with the decomposition of $R_{k+1}$, we use \cref{SchauderBoundLemma} and get 
\[
||\tilde{\phi}||_{C^{2,\alpha}_{\eps}(M)} = O(\eps^{k+1})
\] 
This finishes the induction. \qed \nl \nl
\noindent As a result, we have the following corollaries for $k = 1$
\begin{corollary} \label{C3AlphaExpansion}
For $u_{\eps}$ a solution to Allen--Cahn with Dirichlet data on $Y = \partial M$ a $C^{4,\alpha}$ hypersurface, we have that 
\begin{align*}
u_{\eps}(s,t) &= g_{\eps}(t) + \eps H_Y(s) w_{\eps}(t) + \phi(s,t) \\
||\phi||_{C^{2,\alpha}_{\eps}(M)} &= O(\eps^2)
\end{align*}
\end{corollary}
\noindent Similarly for $k = 2$, we have
\begin{corollary} \label{C5AlphaExpansion}
For $u_{\eps}$ a solution to Allen--Cahn with Dirichlet data on $Y = \partial M$ a $C^{5,\alpha}$ hypersurface, we have that 
\begin{align*}
u_{\eps}^+(s,t) &= g_{\eps}(t) + \eps H_Y(s) w_{\eps}(t) + \eps^2 \left[\dot{H}_0(s) \tau_{\eps}(t) + H_Y^2 \rho_{\eps}(t) + \frac{1}{2} H_Y^2(s) \kappa_{\eps}(t) \right] + \phi \\
||\phi||_{C^{2,\alpha}_{\eps}(M)} &= O(\eps^3)
\end{align*}
\end{corollary}
\section{Proof of theorem \ref{ProjectionTheorem}} \label{ProofOfProjectionSection}
In this section, we work in the closed setting. Consider $Y$ a minimal hypersurface, and perturbations $\eta: Y \to \R$, with $Y_{\eta}$ defined as in section \ref{Results}. Recall the definition of $u_{\eps, \eta}$ \eqref{PastedSolution}:
\[
u_{\eps, \eta} = \begin{cases}
u_{\eps, \eta}^+(p) & p \in M^+ \\
-u_{\eps, \eta}^-(p) & p \in M^-
\end{cases}
\]
We aim to prove \cref{ProjectionTheorem}: 
\ProjectionTheorem*
\noindent As a corollary, we can describe the horizontal variation of the solutions constructed in Pacard--Ritore ([\cite{pacard2003constant}, Thm 4.1], [\cite{pacard2012role}, Thm 1.1]). We recall their notation:
\begin{align*}
u_{\eps}(t) &:= \tanh \left( \frac{t}{\sqrt{2}}\right) \\
\overline{u}(y,t) &= u_{\eps}(t - \zeta(y)) + v(y,t) 
\end{align*}
in our notation $v(y,t) \leftrightarrow \phi(s,t)$ and $\zeta(y) \leftrightarrow \eta(s)$. With this, we have 
\begin{corollary} \label{PacardRitoreSolutionCharacterization}
Suppose $\zeta$ is the perturbation constructed in [\cite{pacard2012role}, Thm 3.33] and $v$ the solution to \eqref{ACEquation} with $v^{-1}(0) = Y_{\zeta}$ and $Y$ non-degenerate, minimal, and separating. Then 
\[
\int_{\R} v(s,t) \dot{\bg}_{\eps}(t) dt = \frac{\sqrt{2}}{3} \zeta(s) + O(\eps^{1 + 2\beta})
\]
with error in $C^{2,\alpha}(Y)$. In particular
\[
\int_{\R} \Delta_Y v(s,t) \dot{\bg}_{\eps}(t) = \frac{\sqrt{2}}{3} J_Y(\zeta) + O(\eps^{1 + 2\beta})
\]
\end{corollary}
\noindent \rmk 
\begin{itemize}
\item Note that 
\[
J_Y = \Delta_Y + (|A_Y|^2 + \Ric_g(\nu, \nu))
\]
and we show that 
\[
\int_{\R} (|A_Y|^2 + \Ric_g(\nu, \nu)) v \dot{\bg}_{\eps} = o(\eps^{1 + 2 \beta})
\]
in $C^{\alpha}(Y)$. Thus we could replace $\Delta_Y v$ with $J_Y v$ on the left hand side of \ref{PacardRitoreSolutionCharacterization}

\item The corollary tells us that the Pacard--Ritore solutions \textbf{have horizontal variation as large as the perturbation}, $\zeta(s)$, off of the initial $Y$ minimal. 
\end{itemize}
The proof is essentially the same as \ref{ImprovedNeumannMinimal}, but we have to confront the low regularity of $Y_{\eta}$ given that $\eta \in C^{2,\alpha}(Y)$. We do this by pulling back $u_{\eps, \eta}^+$  to $Y \times [0, -\omega\eps \ln(\eps))$ and then showing:
\begin{equation} \label{EtaNormalDerivative}
\nu^+ (u_{\eps, \eta}^+) = \frac{1}{\eps \sqrt{2}} + \kappa_0 \eps [\Ric(\nu, \nu) + |A_Y|^2] + \left[\sigma_0 J_Y(\eta) + \int_0^{-\omega\eps \ln(\eps)} \Delta_Y(u_{\eps, \eta}) \dot{\bg}_{\eps}\right] + \tilde{O}(\eta^2, \eps^2)
\end{equation}
\subsection{Set up}
For $Y_{\eta}$ as in \eqref{EtaPerturbationEquation}, we consider the decomposition of $M = M_{\eta}^+ \cup_{Y_{\eta}} M_{\eta}^-$ and $u_{\eps, \eta}^{\pm}$ the minimizers on $M_{\eta}^{\pm}$. With Fermi coordinates about $Y$, define
\begin{align*}
\Phi&: Y \times (\omega\eps \ln(\eps)^2, - \omega \eps \ln(\eps)^2) \to Y \times ( \omega\eps \ln(\eps)^2, - \omega \eps \ln(\eps)^2) \\
\Phi(s,t) &:= F\left( s, t + \eta(s) \zeta\left(\frac{-4t}{\omega\eps \ln(\eps)^2} \right) \right)
\end{align*}
for $F$ as in \eqref{FMapGeneral} and where $\zeta(t)$ is the standard bump function which is $1$ on $(-1,1)$ and goes to zero outside of $[-2,2]$. Note the factor of $\ln(\eps)^2$ so that $\Phi$ restricted to $Y \times (\omega \eps \ln(\eps), - \omega \eps \ln(\eps))$ is a diffeomorphism onto its image.
\begin{figure}[h!]
\centering
\includegraphics[scale=0.3]{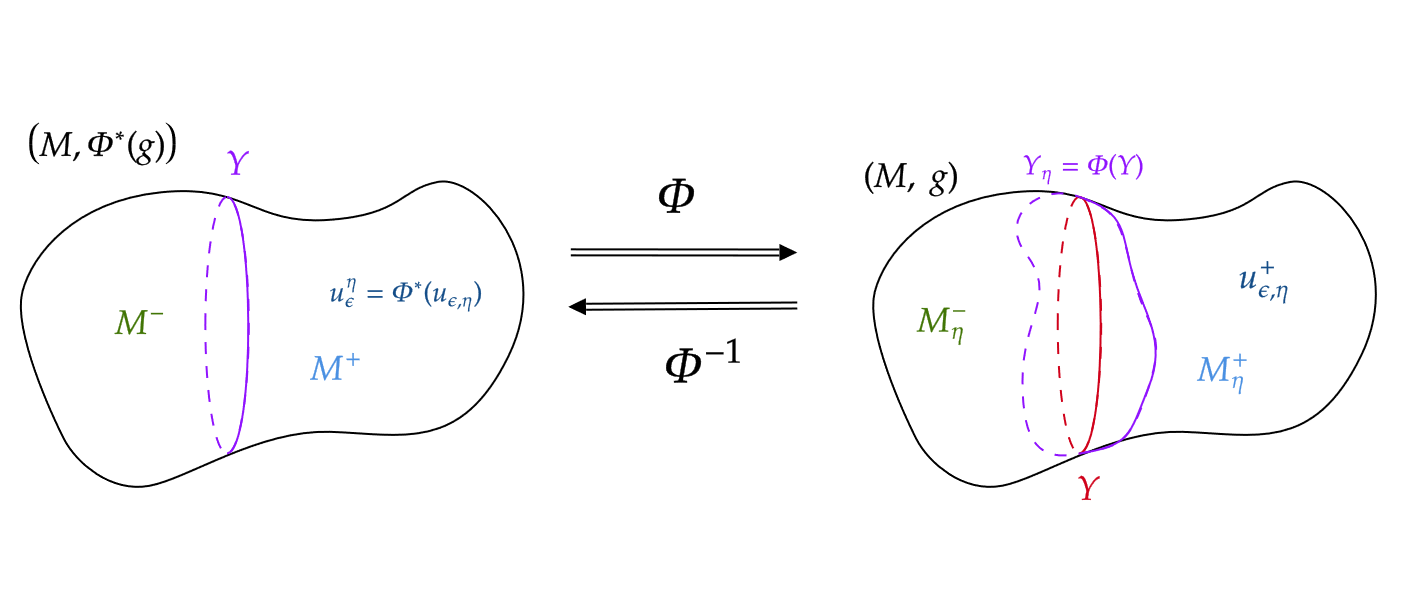}
\caption{$\Phi$ map describing our set up}
\label{fig:acphisetup}
\end{figure}
\nl\noindent
In fact, on this subdomain
\[
|t| < - \omega\eps \ln(\eps) \implies \Phi(s,t) = (s, t + \eta(s))
\]
We pull back $u_{\eps, \eta}^{+}$ by this function and compute the Allen--Cahn equation under this pullback. We define (dropping the $\pm$ notation)
\begin{align*}
u_{\eps}^{\eta}&: Y \times [0, - \omega\eps \ln(\eps)) \to \R \\
&:= \Phi^*(u_{\eps, \eta})(s,t)
\end{align*}
so that 
\begin{align*}
\Phi^*(W'(u_{\eps, \eta})) &= W'(u_{\eps}^{\eta}) \\
u_{\eps}^{\eta}(s,0) &\equiv 0
\end{align*}
and 
\[
\Phi^*(\Delta_g u_{\eps, \eta}^+) = \Delta_{\Phi^*(g)} u_{\eps}^{\eta}
\]
using diffeomorphism invariance of the laplacian. Instead of computing $\Phi^*(g)$ in coordinates, we push forward $\Delta_g$ to $(M, \Phi^*(g))$, as a differential operator, by $\Phi^{-1}$, i.e.
\[
\Delta_{\Phi^*(g)} = (\Phi^{-1})_*(\Delta_g)
\]
We first expand $\Delta_g$ on $(M, g)$, recalling \eqref{LaplacianDecomposition}
\begin{align*}
\Delta_g &= \Delta_t - H_t \partial_t + \partial_t^2 
\end{align*}
Where $\Delta_{t_0}$ denotes the laplacian on $Y_{t_0} := Y + t_0 \nu$, i.e. the set of a signed distance $t_0$ from $Y$. We now compute $\Delta_{\Phi^*(g)} = (\Phi^{-1})_*(\Delta_g)$ by pushing forward each summand:
\begin{align*}
(\Phi^{-1})_*(H_t \partial_t)&= H_{t + \eta} \partial_t \\
(\Phi^{-1})_*(\partial_t^2) &= \partial_t^2 \\
(\Phi^{-1})_*(\Delta_t) &= \Delta_{t + \eta} + E_{\eta}  \\
\implies (\Phi^{-1})_*(\Delta_g) \Big|_t &= \Delta_g \Big|_{t + \eta} + E_{\eta}
\end{align*}
where in the last line, $\Delta_g \Big|_{t + \eta}$ denotes the ambient laplacian on $M$ but with metric coefficients evaluated at the point $(s, t + \eta(s))$. We also define
\begin{align} \label{EDefinition} \nonumber
E_{\eta}&:= g^{ij}(s,t + \eta(s))\left[-\eta_i(s) \partial_t \partial_{s_j} - \eta_j(s) \partial_t \partial_{s_i} - \eta_{ij}(s) \partial_t + \eta_i(s) \eta_j(s) \partial_t^2\right] \\
&\;= - \Delta_{t + \eta}(\eta) \partial_t - 2 \nabla^{t + \eta}(\eta) \partial_t + |\nabla^{t + \eta} \eta|^2 \partial_t^2
\end{align}
From hereon, we only consider $u_{\eps}^{\eta}$ restricted to $Y \times [0, - \omega\eps \ln(\eps))$, and we rewrite the pulled back allen-cahn equation as:
\[
\eps^2 [\Delta_{t + \eta} - H_{t + \eta} \partial_t + \partial_t^2 + E_{\eta}](u^{\eta}_{\eps}) = W'(u_{\eps}^{\eta})
\]
Now we decompose (using that $Y$ is minimal and inspired by \cref{C5AlphaExpansion})
\begin{equation} \label{uEtaInitialDecomposition}
u_{\eps}^{\eta} = \bg_{\eps}(t) + \eps^2 \dot{H}_0(s) \overline{\tau}_{\eps}(t) + \phi^{\eta}(s,t)
\end{equation}
where 
\[
\ddot{\tau}(t) - W''(g(t)) \tau(t) = t \dot{g}(t)
\]
on $\R^+$. From hereon, we label $\phi^{\eta} =: \phi$, and our pulled back Allen--Cahn equation becomes
\begin{align} \label{PulledBackAllenCahn}
L_{\eps, t + \eta}(\phi) &= \eps [\Delta_{t + \eta}(\eta) + \dot{H}_0 \eta] \dot{\bg}_{\eps} + |\nabla^{t + \eta} \eta|^2 \ddot{\bg}_{\eps}  \\ \nonumber
&+ \eps^3 \left( \frac{H_{t + \eta}(s) - \dot{H}_0(s) (t + \eta)}{(t + \eta)^2} \right) \left(\frac{t + \eta}{\eps}\right)^2 \dot{\bg}_{\eps} \\ \nonumber
& -\eps^4 \Delta_{t + \eta}(\dot{H}_0) \overline{\tau}_{\eps} + \eps^3 H_{t + \eta} \dot{H}_0 \dot{\overline{\tau}}_{\eps} + \eps^3 \Delta_{t + \eta}(\eta) \dot{H}_0 \dot{\overline{\tau}}_{\eps} + 2 \eps^3 \nabla^{t+\eta}(\eta)(\dot{H}_0) \dot{\overline{\tau}}_{\eps} \\ \nonumber
& - \eps^2 |\nabla^{t + \eta} \eta|^2 \dot{H}_0 \ddot{\overline{\tau}}_{\eps} + \frac{1}{2} W'''(\bg_{\eps}) \eps^4 \dot{H}_0^2 \overline{\tau}_{\eps}^2 + \eps^6 \dot{H}_0^3 \overline{\tau}_{\eps}^3 \\ \nonumber
& + R(\phi) + O(\eps^{\omega})
\end{align}
where
\begin{align} \nonumber
L_{\eps, r} &:= \eps^2 (\Delta_r - H_r \partial_t + \partial_t^2 - W''(\bg_{\eps}(t))) \\ \label{RError}
R(\phi) &:= \eps^2 E_{\eta}(\phi) - F_0(\phi) \\ \label{F0Error}
F_0(\phi) &:= W'''(\bg_{\eps}) \eps^2 \dot{H}_0 \overline{\tau}_{\eps} \phi + \left[3 \eps^4 \dot{H}_0^2 \tau_{\eps}^2 + \frac{1}{2} W'''(\bg_{\eps}) \right]\phi^2 + \phi^3
\end{align}
where $F_0(\phi)$ is the error term from expanding $W'(u)$. Here, all of the $O(\eps^{\omega})$ terms come from replacing $g_{\eps} \to \bg_{\eps}$ and the like. We abbreviate the right hand side of equation \eqref{PulledBackAllenCahn} as $G(\phi)$. Note that we have (and will continue to) abuse notation with $L_{\eps, t} = L_{\eps}$.

\subsection{Estimates on $\phi$} \label{PhiEtaEstimates}
Again using \cref{SchauderBoundLemma} and equation \eqref{decayAway}, we have
\begin{align*}
||\phi||_{C^{2,\alpha}_{\eps}(M)} &\leq K ||L_{\eps} \phi||_{C^{\alpha}_{\eps}(M)} \\
& \leq \left(\frac{2}{\delta}\right)^{\alpha} [||L_{\eps} \phi||_{C^{\alpha}_{\eps}(t < \delta)} + ||L_{\eps} \phi||_{C^{\alpha}_{\eps}(t > \delta/2)}] \\
& \leq K ||L_{\eps} \phi||_{C^{\alpha}_{\eps}(t < \delta)} + O(\eps^{\omega}) \\
& \leq K ||L_{\eps, t + \eta} \phi||_{C^{\alpha}_{\eps}(t < \delta)} + ||\eps^2 (\Delta_{t + \eta} - \Delta_t) \phi||_{C^{\alpha}_{\eps}(t < \delta)} + ||\eps^2 (H_{t + \eta} - H_t) \phi_t||_{C^{\alpha}_{\eps}(t < \delta)} + O(\eps^{\omega})
\end{align*}
having used \eqref{decayAway} to bound $||L_{\eps} \phi||_{C^{\alpha}_{\eps}(t > \delta/2)}$. With $||\eta||_{C^{2,\alpha}(Y)} \leq K \eps^{1 + \beta}$, this implies
\begin{align} \nonumber
||\phi||_{C^{2,\alpha}_{\eps}(M)} & \leq K ||L_{\eps, t + \eta} \phi||_{C^{\alpha}_{\eps}(t < \delta)} + \eps^{1 + \beta} ||\phi||_{C^{2,\alpha}_{\eps}(M)} + O(\eps^{\omega}) \\
\implies ||\phi||_{C^{2,\alpha}_{\eps}(M)} & \leq K \left[ ||\eps [\Delta_{t + \eta} + \dot{H}_0](\eta) \dot{\bg}_{\eps}||_{C^{\alpha}_{\eps}(t < \delta)} +  || \; |\nabla^{t + \eta} \eta|^2 \ddot{\bg}_{\eps}||_{C^{\alpha}_{\eps}(t > \delta)} \right] + O(\eps^{3}) \\ \label{InitialPhiEtaEstimate}
& \leq O( \eps^{2 + \beta})
\end{align}
\subsection{Proof of \cref{ProjectionTheorem}}
Now as in \S \ref{ProofOfNeumannResult}, we again decompose $L_{\eps, t + \eta}$, multiply by $\dot{\bg}_{\eps}$, integrate, and extract the normal derivative:
\begin{align*}
\int_0^{- \omega \eps \ln(\eps)} L_{\eps, t + \eta}(\phi) \dot{\bg}_{\eps} &= - \eps^2 \sigma \phi_t(s,0) + \eps^2 \int_{0}^{- \omega\eps \ln(\eps)} \Delta_Y(\phi) \dot{\bg}_{\eps} +  O(\eps^{4+\beta-\alpha}) 
\end{align*}
where the above holds in $C^{\alpha}(Y)$, having used \eqref{InitialPhiEtaEstimate}. Similarly
\begin{align*}
\int_0^{- \omega \eps \ln(\eps)} G(\phi) \dot{\bg}_{\eps} &= \eps^2 \dot{H}_0 \eta \sigma_0 + \eps^2 \Delta_0 (\eta) \sigma_0 + O(\eps^{3 + 2 \beta})
\end{align*}
with error terms hold in $C^{\alpha}(Y)$. The details are sketched in the appendix (section \S \ref{IntegratingEtaEquation}). Because $||\eta||_{C^{2,\alpha}(Y)} \leq K \eps^{1 + \beta}$, we see that in terms of order of $\eps$
\[
\phi_t(s,0) = \underbracket{\sqrt{2}\int_0^{- \omega \eps \ln(\eps)} \Delta_Y(\phi) \dot{\bg}_{\eps}(t) dt}_{O(\eps^{1 + \beta - \alpha})} -  \underbracket{\frac{2}{3}[\Delta_0 + \dot{H}_0](\eta)}_{O(\eps^{1 + \beta})} + O(\eps^{1 + 2 \beta})
\]
where the above asymptotics hold in $C^{\alpha}(Y)$. We frame this as
\begin{equation} \label{EtaNormalDerivative}
\Big|\Big|\phi_t(s,0) - \sqrt{2} \int \Delta_Y(\phi) \dot{\bg}_{\eps} + \frac{2}{3} J_Y(\eta)\Big|\Big||_{C^{\alpha}(Y)} =  O(\eps^{1 + 2 \beta})
\end{equation}
We now note that $\partial_t$ is comparable to $(\Phi^{-1})_*(\nu_{\eta})$ (see section \S \ref{EtaNormal}), i.e. the normal vector for $Y$ and that of $Y_{\eta}$ (translated to $Y$) are comparable since $\eta$ is small:
\begin{align*}
\nu_{\eta} &= (1 + A(s))\partial_t + B^i(s) \partial_{s_i} \Big|_{t = \eta(s)}\\
||A(s)||_{C^{\alpha}(Y)} &\leq C ||\eta||_{C^{1,\alpha}(Y)}^2 \\
||B^i(s)||_{C^{\alpha}(Y)}& \leq C ||\eta||_{C^{1,\alpha}(Y)} \\
\implies (\Phi^{-1})_*(\nu_{\eta}) &= (1 + \tilde{A}(s)) \partial_t + \tilde{B}^i(s) \partial_{s_i} \Big|_{t = 0} \\
||\tilde{A}(s)||_{C^{\alpha}(Y)} &\leq C ||\eta||_{C^{1,\alpha}(Y)}^2 \\
||\tilde{B}^i(s)||_{C^{\alpha}(Y)}& \leq C ||\eta||_{C^{1,\alpha}(Y)} 
\end{align*}
as such 
\[
||(\Phi^{-1})_*(\nu_{\eta})(\phi) - \phi_t(s,0)||_{C^{\alpha}(Y)} = O(\eps^{ 3 + 3 \beta - \alpha})
\]
(recall that $\phi_{s_i}(s,0) \equiv 0$) and so 
\begin{align*}
(\Phi^{-1})_*(\nu_{\eta})(\phi) &= \phi_t(s,0) + O(\eps^{3 + 3 \beta - \alpha}) \\ 
& = \sqrt{2} \int_0^{- \omega \eps \ln(\eps)} \Delta_Y(\phi) \dot{\bg}_{\eps}(t) dt - \frac{2}{3} J_Y(\eta) + O(\eps^{1 + 2 \beta})
\end{align*}
To prove theorem \ref{ProjectionTheorem}, we note that if $u_{\eps, \eta}$ is $C^1$ across $Y_{\eta} = u_{\eps, \eta}^{-1}(0)$, then the Neumann data match. If we take 
\begin{align*}
u_{\eps, \eta} &: Y \times ( \omega \eps \ln(\eps)^2, - \omega \eps \ln(\eps)^2) \to \R \\
\overline{u}_{\eps}^{\eta} &:= \Phi^*(u_{\eps, \eta}): Y \times ( \omega \eps \ln(\eps), - \omega \eps \ln(\eps)) \to \R \\
\overline{u}_{\eps}^{\eta} &= \bg_{\eps}(t) + \eps^2 \dot{H}_0(s) \overline{\tau}_{\eps}(t) + \overline{\phi}(s,t) \\
\overline{\tau}_{\eps}(t) &:= \begin{cases}
\tau_{\eps}(t) & t \geq 0 \\
- \tau_{\eps}(-t) & t < 0
\end{cases} \\
\overline{\phi}(s,t)&:= \begin{cases}
\phi^{+}(s,t) & t \geq 0 \\
\phi^{-}(s,t) & t < 0
\end{cases}
\end{align*}
where $\phi^{\pm}$ are the same functions as in \eqref{uEtaInitialDecomposition} with the $\pm$ made explicit to represent working on $M^{\pm}$ (i.e. $t > 0$ or $t < 0$). With the above, $\nu_{\eta}(u_{\eps, \eta}) = (\Phi_{\eta})^{-1}_*(\nu_{\eta}) (u_{\eps}^{\eta})$ is well defined, and 
\begin{align*}
0 &= \nu_{\eta}(u_{\eps, \eta})\Big|_{t = \eta(s)^+} - \nu_{\eta}(u_{\eps, \eta}) \Big|_{t = \eta(s)^-} \\
\iff 0 &= (\Phi_{\eta})^{-1}_*(\nu_{\eta}) (u_{\eps}^{\eta})\Big|_{t = 0^+} - (\Phi_{\eta})^{-1}_*(\nu_{\eta}) (u_{\eps}^{\eta}) \Big|_{t = 0^-} \\
& = -\frac{4}{3}J_Y(\eta) + \sqrt{2}\int_{\omega\eps \ln(\eps)}^{- \omega \eps \ln(\eps)} \Delta_Y(\overline{\phi}) \dot{\bg}_{\eps}(t) dt + O(\eps^{1 + 2 \beta}) \\
\implies \int_{\omega\eps \ln(\eps)}^{- \omega \eps \ln(\eps)} \Delta_Y(\overline{\phi}) \dot{\bg}_{\eps}(t) dt &= \frac{2 \sqrt{2}}{3}\sigma_0 J_Y(\eta) + \tilde{O}(\eps^{1 + 2 \beta})
\end{align*}
where the error holds in $C^{\alpha}(Y)$.
%\implies \int_{\omega \eps \ln(\eps)}^{- \omega \eps \ln(\eps)} u_{\eps, \eta}(s,t) \dot{g}_{\eps}(t) dt &= \sigma_0 \eta(s) + \tilde{O}(\eps^{1 + 2 \beta})
Now note that
\[
\Big| \Big|\int_{ \omega\eps \ln(\eps)}^{- \omega\eps \ln(\eps)} \dot{H}_0 \overline{\phi} \dot{\bg}_{\eps} dt\Big| \Big|_{C^{\alpha}(Y)} = O(\eps^{3 + \beta - \alpha} )
\]
so that we can write the above as 
\[
J_Y\left(\int_{\omega\eps \ln(\eps)}^{- \omega\eps \ln(\eps)} \overline{\phi} \dot{\bg}_{\eps}(t) dt\right) = 2 \sigma_0 J_Y(\eta) + O(\eps^{1 + 2 \beta})
\]
with error in $C^{\alpha}(Y)$. We now substitute $\int u_{\eps, \eta} \dot{\bg}_{\eps}$ for $\int \phi \dot{\bg}_{\eps}$ at the cost of negligible error since
\begin{align*}
\Big|\int_{ \omega\eps \ln(\eps)}^{- \omega \eps \ln(\eps)} g_{\eps}(t) \dot{g}_{\eps}(t) dt \Big| &= O(\eps^{\omega}) \\
\Big| \Big|\int_{ \omega\eps \ln(\eps)}^{- \omega \eps \ln(\eps)} \eps^2 \dot{H}_0(s) \tau_{\eps}(t) \dot{g}_{\eps}(t) \Big| \Big|_{C^{\alpha}(Y)} &= O(\eps^3)
\end{align*}
and the same order of bound holds if we replace $g_{\eps} \to \bg_{\eps}, \tau_{\eps} \to \overline{\tau}_{\eps}$. Furthermore, if $Y$ is non-degenerate, we invert both sides by $J_Y$
\[
\int_{ \omega \eps \ln(\eps)}^{- \omega\eps \ln(\eps)} u_{\eps, \eta}(s,t) \dot{\bg}_{\eps}(t) dt = 2\sigma_0 \eta(s) + O(\eps^{1 + 2 \beta})
\]
where the error holds in $C^{2,\alpha}(Y)$. This concludes the theorem. \qed

\section{Future Work}
The following generalizations and next steps are of interest:
\begin{itemize}
\item Can we reproduce the multiplicity one version of Wang--Wei [\cite{wang2019second}, Thm 1.1] using these techniques? More specifically, instead of perturbing the $0$ set of Allen--Cahn solutions by some function $h(s)$, can we bootstrap a system of equations using Schauder estimates from \cref{SchauderBoundLemma} to get regularity of $H_{Y_{\eta}}$? This would require showing a better bound than $||\eta||_{C^{2,\alpha}(Y)} \leq K \eps^{1 + \beta}$ when $u_{\eps, \eta}$ being a solution to \eqref{ACEquation}

\item Can we show that 
\[
\int u_{\eps, \eta} \dot{g}_{\eps} = \sigma_0 \eta(s) + \tilde{O}(\eps^2, \eta^2)
\]
holds for all such perturbations, $\eta$, off of $Y$ minimal? I.e. not \textit{just} when $\eta$ gives rise to a 2-sided solution to \eqref{ACEquation} via energy minimization.

\begin{itemize}
\item  Given the above equality, can one reprove the result of Pacard--Ritore using just energy minimization and the local perturbation $\eta$? We note the recent work of De Phillipis and Pigati [\cite{de2022non}, Thm 5] who prove the result using purely energy and gradient flow methods. We ask if an alternate solution can be performed by better understanding the Dirichlet-to-Neumann map for arbitrary perturbation.

\item In particular, consider $Y$ a non-degenerate minimal hypersurface and the maps
\begin{align*}
    \mathcal{P}(\eta) &= \frac{\partial u^+_{\eps, \eta}}{\partial \nu_{\eta}} - \frac{\partial u^-_{\eps, \eta}}{\partial \nu_{\eta}} \\
    F(\eta) &:= - \frac{1}{4 \sigma_0}J_Y^{-1}[ \mathcal{P}(\eta) - 4\sigma_0 J_Y(\eta)] = -J_Y^{-1}(\tilde{R}_{Y}(\eta))
\end{align*}
with the hopes of showing that $F: C^{2,\alpha}(Y) \to C^{2,\alpha}(Y)$ is a contraction.

\item The above boils down to showing that 
\[
\Big| \Big| \int_0^{- \omega \eps \ln(\eps)} \Delta_Y(\phi) \dot{g}_{\eps}dt\Big| \Big|_{C^{\alpha}(Y)} = o(\eps^{1 + \beta})
\]
\end{itemize}
\item During the course of this paper, the author conjectured the following: 
\begin{conjecture}\label{HalfSpaceCharacterization}
Let $\phi: \R^{n-1} \times \R^+ \to \R$. Suppose that $\phi(s,0) \equiv 0$. Let $f: \R^{n-1} \to \R$, and suppose 
\[
[\Delta_{\R^n} - W''(g)](\phi) = f(s) \dot{g}(t)
\]
then there exists $c \in \R$ such that $f(s) \equiv c \in \R$ and $\phi(s,t) = c w(t)$ for $w(t)$ as in \eqref{wEquation}.
\end{conjecture}
\end{itemize}

\section{Appendix}
\subsection{Lemma on $L^* = \Delta_{\R^n} + \partial_t^2 - W''(g)$ for $\R^n \times \R^+$}
\begin{lemma} \label{HalfSpace}
For $\phi \in C^{1}(\R^n \times \R^+)$, suppose $\phi(s,0) \equiv 0$ and $L^*(\phi) = 0$ on $\R^n \times \R^+$. Then $\phi \equiv 0$
\end{lemma}
\noindent \Pf The proof is a slight extension of the well known classification of $\ker(L)$ on $\R^n \times \R$. See [\cite{pacard2012role}, Lemma 3.7] for reference. Because of the Dirichlet condition at $t = 0$, consider the odd reflection
\[
\tilde{\phi}(s,t) = \begin{cases}
\phi(s,t) & t \geq 0 \\
-\phi(s,-t) & t < 0
\end{cases}
\]
then $\tilde{\phi}(s,t)$ is a $C^1$ solution to $L^*$. By the maximum principle, $\tilde{\phi}$ converges to $0$ exponentially and uniformly as $t \to \pm \infty$. Thus it is in $L^2$ and via an energy argument (again [\cite{pacard2012role}, Lemma 3.7]), we see that
\[
\tilde{\phi}(s,t) = c \dot{g}(t)
\]
but $\phi(s,0) = 0$ so $c = 0$. \qed 

\subsection{Boundedness of $\tau$ in \cref{SchauderBoundLemma}} \label{TauBounded}
Recall that we have a sequence $\{f_j\}: M \to \R$ and $p_j$ such that $|f_j(p_j)| = ||f_j||_{C^0(M)} = 1$ and $||L_{\eps} f_j||_{C^{\alpha}_{\eps}(M)} \leq j^{-1}$. We want to show that $\text{dist}(Y, p_j) < \kappa \eps_j$ for some $\kappa$ independent of $\eps_j$. \nl \nl
In terms of scaled fermi coordinates, $(\sigma, \tau)$, we have 
\[
L_{\eps} = \Delta_{g_{\eps}} - W''(\bg(\tau))
\]
Consider $\phi = 1$ for which 
\[
L(1) = - W''(\bg) = 1 - 3\bg^2
\]
We see that this $L(1) < -1$ for all $\tau > \text{arctanh}\left(\sqrt{\frac{2}{3}} \right) =: c_0$. Moreover, $1 \geq |\phi_j|$ by the normalization. We now apply the maximum principle to $L$ and $(1 \pm \phi)$ on the open set $U = \{p \; | \; \text{dist}(p, Y) > c_0 \eps \}$. This tells us that for $\tau > c_0$ (i.e. $t > c_0 \eps$), $1 \pm \phi_j$ achieves its minimum on the boundary of $\{\tau > c_0\}$. Immediately, this tells us that we can choose $p_j = (q_j, t_j)$ for some $0 \leq t_j < c_0 \eps$. \nl \nl
Similarly, we can show that $\tau_j = \eps^{-1} t_j \geq \tau_0 > 0$. Recentering $\phi$ at $(q_j, 0)$ and using $(\sigma, \tau)$ coordinates, we have
\begin{align*}
|\partial_\tau \phi_j(\sigma, \tau)| &\leq ||\phi_j||_{C^{2,\alpha}_{\eps}(M))} \\
&\leq K (||L \phi_j||_{C^{\alpha}_{\eps}(M)} + ||\phi_j||_{C^0(M)}) \\
&\leq K (o(1) + 1) \leq 2K
\end{align*}
so that because $\phi_j(\sigma,0) \equiv 0$ for all $\sigma$, we have that 
\[
|\phi_j(\sigma_j, \tau_j)| \geq 1/2 \implies \tau_j \geq \frac{1}{4K}
\]
where $K$ is the Schauder constant and independent of $j$ and $\eps$. This tells us that 
\[
0 < \frac{1}{4K} \leq \tau_j \leq c_0
\]
so there exists a convergent subsequence of $\{\tau_j\}$ which converges to $0 < \tau < \infty$.
\subsection{Normal for $Y_{\eta}$} \label{EtaNormal}
In this section, we show that for $\eta$ a perturbation with $||\eta||_{C^{2,\alpha}}|| \leq K \eps^{1 + \beta}$, we have
\begin{align*}
\nu_{Y_{\eta}} &= a^t(s) \partial_t + a^i(s) \partial_{s_i} \\
||a^t - 1||_{C^{2,\alpha}(Y)} &\leq C \eps^{2 + 2\beta} \\
||a^i ||_{C^{2,\alpha}(Y)} &\leq C \eps^{1 + \beta}
\end{align*}
\begin{lemma}
For any $\eta \in C^{2,\alpha}(Y)$ and $||\eta||_{C^{2,\alpha}(Y)} \leq K \eps^{1 + \beta}$, there exists $C > 0$ so that, the normal derivative to $Y_{\eta}$ expands as 
\begin{align*}
\nu_{\eta} &= a^t(s) \partial_t + a^i(s) \partial_{s_i}  \\
||a^t(s) - 1||_{C^{1,\alpha}(Y)} &\leq C ||\eta||_{C^{2,\alpha}(Y)}^2 \\
||a^i(s)||_{C^{1,\alpha}(Y)} &\leq C ||\eta||_{C^{2,\alpha}(Y)}
\end{align*}
\end{lemma}
\noindent \Pf In coordinates, we compute the tangent basis for $Y_{\eta}$ as 
\[
v_i = \partial_{s_i} + \eta_i \partial_t \Big|_{(s, t = \eta)}
\]
Let $g(\eta)_{ij} := g(v_i, v_j)$ and $g(\eta)^{ij}$ be the corresponding inverse. Then
\begin{align*}
g(\eta)_{ij} &= g_{ij} \Big|_{(s, t = \eta)} + \eta_i \eta_j  \\
&= \delta_{ij} + \eta A_{ij} + \eta_i \eta_j \\
&= \delta_{ij} + \eta A_{ij} + O((D \eta)^2)
\end{align*}
so that for
\begin{align*}
w_{\eta} &:= \partial_t - \Pi_{TY}(\partial_t) \\
& = \partial_t - g^{ij} g(\partial_t, v_i) v_j \\
& = \partial_t - (\delta_{ij} - \eta A_{ij} + \tilde{O}((D \eta)^2)) \eta_i (\partial_{s_j} + \eta_j \partial_t) \\
& = (1 + \tilde{O}(\eta \eta_i \eta_j, (D \eta)^4) \partial_t - ( \eta_i \delta_{ij} + \tilde{O}( \eta D \eta, (D \eta)^3)) \partial_j 
\end{align*}
The third line comes from expanding the metric in fermi coordinates over $Y$ and evaluating at $t = \eta$. We compute
\[
||w_{\eta}||^2 = 1 + \tilde{O}((D \eta)^2) \implies ||w_{\eta}||^{-1} = 1 + \tilde{O}((D\eta)^2)
\]
and so 
\begin{align*}
\nu_{\eta} &= \frac{w_{\eta}}{||w_{\eta}||} = a^t(s) \partial_t + B^i(s) \partial_{s_i} \\
||a^t - 1||_{C^{1,\alpha}(Y)} &\leq C || \eta||_{C^{2,\alpha}(Y)}^2 \\
||B^i ||_{C^{1,\alpha}(Y)} &\leq C || \eta||_{C^{2,\alpha}(Y)}
\end{align*}
\qed \nl 
\noindent We now note that for the diffeomorphism
\begin{align*}
\Phi&: Y \times (\omega \eps \ln(\eps), - \omega \eps \ln(\eps)) \to Y \times (-\omega \eps \ln(\eps)^2, \omega \eps \ln(\eps)^2) \\
\Phi(s,t) &= (s, t + \eta) \\
\implies (\Phi^{-1})_*(\nu_{\eta}) &= (1 + A(s)) (\Phi^{-1})_*(\partial_t) + B^i(s) (\Phi^{-1})_*(\partial_{s_i}) \\
&= (1 + A(s)) \partial_t + B^i(s) (\partial_{s_i} + \eta_i \partial_t) \\
&= (1 + A(s) + B^i \eta_i) \partial_t + B^i(s) \partial_{s_i} \\
&= (1 + \tilde{A}(s)) \partial_t + B^i(s) \partial_{s_i} \\
||\tilde{A}(s)||_{C^{1,\alpha}(Y)} & \leq C ||\eta||_{C^{2,\alpha}(Y)}^3 \leq C ||\eta||_{C^{2,\alpha}(Y)}^2
\end{align*}
using that $||\eta||_{C^{2,\alpha}(Y)} = O(\eps^{1 + \beta})$.
% \nl \nl
%
\subsection{Integrating $\eta$ equation for normal derivative } \label{IntegratingEtaEquation}
In this section, we keep track of all the terms in equation \eqref{PulledBackAllenCahn} when integrating against $\dot{\bg}_{\eps}$ to extract $\phi^{\eta}_t(s,0)$. 
\begin{lemma}
$\phi_t^{\eta}(s,0)$ decomposes as 
\[
\phi_t^{\eta}(s,0) = \sqrt{2}\int_0^{- \omega \eps \ln(\eps)} \Delta_Y(\phi) \dot{\bg}_{\eps}(t) dt - \frac{2}{3}J_Y(\eta) + O(\eps^{1 + 2 \beta})
\]
where the error bound holds in $C^{\alpha}(Y)$.
\end{lemma}
\noindent \Pf Recall from \eqref{PulledBackAllenCahn} that we have 
\begin{align} \label{LinEtaEq}
L_{\eps, t + \eta}(\phi) &= \eps [\Delta_{t + \eta}(\eta) + \dot{H}_0 \eta] \dot{\bg}_{\eps} + |\nabla^{t + \eta} \eta|^2 \ddot{\bg}_{\eps}  \\ \nonumber
&+ \eps^3 \left( \frac{H_{t + \eta}(s) - \dot{H}_0(s) (t + \eta)}{(t + \eta)^2} \right) \left(\frac{t + \eta}{\eps}\right)^2 \dot{\bg}_{\eps} \\ \nonumber
& -\eps^4 \Delta_{t + \eta}(\dot{H}_0) \overline{\tau}_{\eps} + \eps^3 H_{t + \eta} \dot{H}_0 \dot{\overline{\tau}}_{\eps} + \eps^3 \Delta_{t + \eta}(\eta) \dot{H}_0 \dot{\overline{\tau}}_{\eps} + 2 \eps^3 \nabla^{t+\eta}(\eta)(\dot{H}_0) \dot{\overline{\tau}}_{\eps} \\ \nonumber
& - \eps^2 |\nabla^{t + \eta} \eta|^2 \dot{H}_0 \ddot{\overline{\tau}}_{\eps} + \frac{1}{2} W'''(\bg_{\eps}) \eps^4 \dot{H}_0^2 \overline{\tau}_{\eps}^2 + \eps^6 \dot{H}_0^3 \overline{\tau}_{\eps}^3 \\ \nonumber
& + R(\phi) + O(\eps^{\omega})\\
& =: G(\phi)
\end{align}
Starting with the left hand side, we first recall
\[
L_{\eps, t + \eta}(\phi) = \eps^2(\Delta_{t + \eta}(\phi^{\eta}) - H_{t + \eta} \phi_t^{\eta} + \phi_{tt}^{\eta}) - W''(\bg_{\eps})\phi^{\eta} 
\] 
we multiply by $\dot{\bg}_{\eps}$ and integrate from $t = 0 \to t = -\omega \eps \ln(\eps)$. From hereon, all integrals will be from $[0, -\omega \eps \ln(\eps))$
\begin{align*}
\int_0^{-\omega\eps \ln(\eps)} L_{\eps, t + \eta}(\phi) \dot{\bg}_{\eps}&= - \eps^2 \sigma\phi^{\eta}_t(s,0) + \eps^2 \int \Delta_{t + \eta}(\phi^{\eta})\dot{\bg}_{\eps} - \eps^2\int H_{t + \eta} \phi_t^{\eta} \dot{\bg}_{\eps}\\
&= - \eps^2 \sigma \phi^{\eta}_t(s,0) + \eps^2 \int \Delta_Y(\phi) \dot{\bg}_{\eps} + O(\eps^{4-\alpha})
\end{align*}
here we've used  \eqref{InitialPhiEtaEstimate} i.e. 
\[
||\phi||_{C^{2,\alpha}_{\eps}(M))} = O(\eps^{2 + \beta})
\]
and $||\int \dot{\bg}_{\eps} H_{t + \eta}||_{C^{\alpha}(Y)} = O(\eps^{1 + \beta})$ since $Y$ is minimal and $||\eta||_{C^{2,\alpha}} = O(\eps^{1 + \beta})$. On the right hand side of \eqref{LinEtaEq}, we have 
\begin{align*}
\int_0^{-\omega\eps \ln(\eps)} G(\phi^{\eta}) \dot{\bg}_{\eps} dt	&=  \int \eps [\Delta_{t + \eta}(\eta) + \dot{H}_0 \eta] \dot{\bg}_{\eps}^2 dt +  \int |\nabla^{t + \eta} \eta|^2 \dot{\bg}_{\eps} \ddot{\bg}_{\eps} dt  \\ \nonumber
&+ \eps^3 \int \left( \frac{H_{t + \eta}(s) - \dot{H}_0(s) (t + \eta)}{(t + \eta)^2} \right) \left(\frac{t + \eta}{\eps}\right)^2 \dot{\bg}_{\eps}^2 dt \\ \nonumber
& -\eps^4 \int \Delta_{t + \eta}(\dot{H}_0) \overline{\tau}_{\eps} \dot{\bg}_{\eps} dt + \eps^3 \int H_{t + \eta} \dot{H}_0 \dot{\overline{\tau}}_{\eps} \dot{\bg}_{\eps} dt + \eps^3 \int \Delta_{t + \eta}(\eta) \dot{H}_0 \dot{\overline{\tau}}_{\eps} \dot{\bg}_{\eps} + 2 \eps^3 \int \nabla^{t+\eta}(\eta)(\dot{H}_0) \dot{\overline{\tau}}_{\eps} \dot{\bg}_{\eps} \\ \nonumber
& - \eps^2 \int |\nabla^{t + \eta} \eta|^2 \dot{H}_0 \ddot{\overline{\tau}}_{\eps} \dot{\bg}_{\eps} + \frac{1}{2} \eps^4 \int W'''(\bg_{\eps}) \dot{H}_0^2 \overline{\tau}_{\eps}^2 \dot{\overline{g}}_{\eps} + \eps^6 \int \dot{H}_0^3 \overline{\tau}_{\eps}^3 \dot{\bg}_{\eps} \\ \nonumber
& + \int R(\phi) \dot{\bg}_{\eps} + O(\eps^{\omega + 1}) \\
\end{align*}
we write this as 
\begin{align*}
\int_0^{-\omega\eps \ln(\eps)} G(\phi^{\eta}) \dot{\bg}_{\eps} dt	&= A_1 + A_2 \\
&+ B_1 \\
&+ C_1 + C_2 + C_3 + C_4 \\
&+ D_1 + D_2 + D_3 \\
&+ E
\end{align*}
With the aim of extracting the leading terms and an appropriate error bounded in $C^{\alpha}(Y)$. We have 
\begin{align*}
A_1 &= \eps \left(\Delta_0(\eta) + \dot{H}_0 \eta \right)\int\dot{\bg}_{\eps}^2 dt + \eps \int [\Delta_{t + \eta} - \Delta_0](\eta) \dot{\bg}_{\eps}^2 dt \\
&= \sigma_0 \eps^2 J_Y(\eta) + \eps \int [\Delta_{t + \eta} - \Delta_0](\eta) \dot{\bg}_{\eps}^2 dt  \\
&= \sigma_0 \eps^2 J_Y(\eta) + O(\eps^3 ||\eta||_{C^{2,\alpha}(Y)}, \eps^2 ||\eta||_{C^{\alpha}(Y)} ||\eta||_{C^{2,\alpha}(Y)}) \\
&= \sigma_0 \eps^2 J_Y(\eta) + O(\eps^{4 + \beta})
\end{align*}
which comes from expanding $\Delta_{t + \eta} - \Delta_0$ in powers of $(t + \eta)$. Similarly
\begin{align*}
A_2 &= \int |\nabla^Y \eta|^2 \dot{\bg}_{\eps} \ddot{\bg}_{\eps} + \int [|\nabla^{t + \eta} \eta|^2 - |\nabla^Y \eta|^2] \dot{\bg}_{\eps} \ddot{\bg}_{\eps} dt \\
&= O(\eps^{3 + 2 \beta})
\end{align*}
which comes from expanding $g^{ij}(s, t + \eta)$ in powers of $(t + \eta)$ and
\begin{align*}
|\nabla^{t + \eta} \eta|^2 &= g^{ij}(s, t + \eta) \eta_i \eta_j \\
&= g^{ij}(s,0) \eta_i \eta_j + [g^{ij}(s, t + \eta) - g^{ij}(s,0)] \eta_i \eta_j
\end{align*}
For $B_1$, we see that
\begin{align*}
||B_1||_{C^0} & \leq \eps^3 \int K \left( \frac{|t| + |\eta|}{\eps} \right) \dot{\bg}_{\eps}^2 dt \\
& \leq O(\eps^4)
\end{align*}
to see the $[\cdot]_{\alpha}$ bound, we write 
\begin{align*}
B_1 &= \eps^3 \int_0^{-\omega \eps \ln(\eps)} \left(\int_0^{t + \eta} [\dot{H}_r(s) - \dot{H}_0(s)] dr\right) \dot{\bg}_{\eps}^2 dt \\
&= \eps^3 \int_0^{-\omega \eps \ln(\eps)} \left(\int_0^{t + \eta} \int_0^r \ddot{H}_w(s) dw dr \right) \dot{\bg}_{\eps}^2 dt \\
\implies [B_1]_{C^{\alpha}(Y)} &= O(\eps^{4})
\end{align*}
For $\{C_i\}$, we compute in a straight forward manner using that $Y \in C^{4,\alpha}$ and satisfies \eqref{UniformGeometry}: 
\begin{align*}
C_1 &= \eps^4 \int \Delta_{t + \eta}(\dot{H}_0) \dot{\overline{\tau}}_{\eps} \dot{\bg}_{\eps} dt \\
&= O(\eps^5) \\
C_2 &= \eps^3 \int H_{t + \eta} \dot{H}_0 \dot{\overline{\tau}}_{\eps} \dot{\bg}_{\eps} dt \\
&= O(\eps^4) \\
C_3 &= \eps^3 \int \Delta_{t + \eta}(\eta) \dot{H}_0 \dot{\overline{\tau}}_{\eps} \dot{\bg}_{\eps} dt\\
&= O(\eps^{5 + \beta}) \\
C_4 &=2 \eps^3 \int \nabla^{t+\eta}(\eta)(\dot{H}_0) \dot{\overline{\tau}}_{\eps} \dot{\bg}_{\eps} dt \\
&= O(\eps^{5 + \beta})
\end{align*}
which is seen from making a change of variables $t \to t/\eps$ to gain another factor of $\eps$, and then noting that the integrals converge and are bounded in $C^{\alpha}(Y)$. \nl \nl
For the $\{D_i\}$ terms, we similarly have:
\begin{align*}
D_1 &= - \eps^2 \int |\nabla^{t + \eta} \eta|^2 \dot{H}_0 \ddot{\overline{\tau}}_{\eps} \dot{\bg}_{\eps} dt \\
&= O(\eps^{5 + 2 \beta}) \\
D_2 &= \frac{1}{2} \eps^4 \int W'''(\bg_{\eps}) \dot{H}_0^2 \overline{\tau}_{\eps}^2 \dot{\bg}_{\eps} dt\\
&= O(\eps^5) \\
D_3 &= \eps^6 \int \dot{H}_0^3 \overline{\tau}_{\eps}^3 \dot{\bg}_{\eps} \\
&= O(\eps^7)
\end{align*}
Finally, recall \eqref{RError} to decompose the $E$ term
\begin{align*}
E_1 &= \int R(\phi) \dot{\bg}_{\eps} \\
&= \int \eps^2 E_{\eta}(\phi) \dot{\bg}_{\eps} - \int F_0(\phi) \dot{\bg}_{\eps}
\end{align*}
And we have 
\begin{align*}
E_{\eta}&= - \Delta_{t + \eta}(\eta) \partial_t - 2 \nabla^{t + \eta}(\eta) \partial_t + |\nabla^{t + \eta} \eta|^2 \partial_t^2 \\
\eps^2 \int (- \Delta_{t + \eta} (\eta)) \phi_t \dot{\bg}_{\eps} & = O(\eps^{5 + 2 \beta - \alpha})\\
- 2 \eps^2 \int \nabla^{t + \eta}(\eta)(\phi_t) \dot{\bg}_{\eps} & = O(\eps^{4 + 2 \beta - \alpha})\\
\eps^2 \int |\nabla^{t + \eta} \eta|^2 \phi_{tt} \dot{\bg}_{\eps}  &= O(\eps^{5 + 3 \beta - \alpha})\\
\implies \eps^2 \int E_{\eta}(\phi) \dot{\bg}_{\eps} &= O(\eps^{4 + 2 \beta - \alpha})
\end{align*}
with bounds holding in $C^{\alpha}(Y)$. Similarly, using the definition of $F_0$ in \eqref{F0Error}
\begin{align*}
F_0(\phi) &:= W'''(\bg_{\eps}) \eps^2 \dot{H}_0 \overline{\tau}_{\eps} \phi + \left[3 \eps^4 \dot{H}_0^2 \overline{\tau}_{\eps}^2 + \frac{1}{2} W'''(\bg_{\eps}) \right]\phi^2 + \phi^3 \\
\eps^2 \int W'''(\bg_{\eps}) \dot{H}_0 \overline{\tau}_{\eps} \phi \dot{\bg}_{\eps} &=  O(\eps^{5 + \beta - \alpha})\\
3 \eps^4 \int \dot{H}_0^2 \overline{\tau}_{\eps}^2 \phi^2 \dot{\bg}_{\eps} &= O(\eps^{9 + 2 \beta - \alpha})\\
\frac{1}{2} \int W'''(\bg_{\eps}) \phi^2 \dot{\bg}_{\eps} &= O(\eps^{5 + 2 \beta - \alpha}) \\
\int \phi^3 \dot{\bg}_{\eps} &= O(\eps^{7 + 3 \beta - \alpha}) \\
\implies \int F_0(\phi) \dot{\bg}_{\eps} &= O(\eps^{5 + \beta - \alpha})
\end{align*}
With this, we've shown that
\[
\int G(\phi^{\eta}) \dot{\bg}_{\eps} dt = \eps^2 \sigma_0 J_Y(\eta) + O(\eps^{3 + 2 \beta})
\]
with error in $C^{\alpha}(Y)$. This finishes the proof. \qed
\subsection{Existence of solutions to $\ddot{F} - W''(g) F = \varphi$} \label{HalfLineSection}
Given $\varphi: [0, \infty) \to \R$ smooth and asymptotically exponentially decaying, consider 
\begin{align} \label{HalfSpaceSystem}
\partial_t^2 F(t) - W''(g(t)) F(t) &= \varphi(t) \\ \nonumber
F(0) &= 0 \\
\lim_{t \to \infty} F(t) &= 0
\end{align}
we reprove the following lemma seen in \cite{mantoulidis2022variational} and proven in [\cite{alikakos1996critical}, Lemma B.1, Remark B.3].
\begin{lemma} \label{HalfSpaceSystemExistence}
Given $\varphi \in C^{\infty}([0, \infty))$ such that 
\[
\exists t_0 > 0, \; K > 0, \; \gamma > 0 \quad \st \quad \forall t > t_0, \qquad |\varphi(t)| \leq K e^{-\gamma t}
\] 
then there exists a smooth solution to the system \eqref{HalfSpaceSystem} with exponential decay.
\end{lemma}
\noindent \Pf Consider the a priori solution of the form 
\[
F(t) = v(t) \dot{g}(t)
\]
where $v(t)$ is to be constructed with $v(0) = 0$. We plug this into \eqref{HalfSpaceSystem}, multiply by $\dot{g}$, and integrate twice to get a general solution of 
\begin{align*}	
v(t) &= b_0 + \int_0^t \dot{g}(s)^{-2} \left[ a_0 + \int_0^s \varphi(r) \dot{g}(r) dr \right] ds
\end{align*}
Using the condition of $v(0) = 0$, we have $b_0 = 0$. Moreover, we can set 
\[
a_0 = - \int_0^{\infty} \varphi(r) \dot{g}(r)
\]
We now show that $v(t)$ is bounded so that $\lim_{t \to \infty} F(t) = 0$. We compute 
\[
\dot{v}(t) = \dot{g}(t)^{-2} \left[ a_0 + \int_0^t \varphi(r) \dot{g}(r) dr \right]
\]
we know that for $t$ large, 
\[
\dot{g}(t)^{-2} \sim e^{\sqrt{2} t}
\]
So it suffices to show that 
\[
|a_0 + \int_0^t \varphi(r) \dot{g}(r)| \leq K e^{-(\beta + \sqrt{2}) t}
\]
for some $K > 0, \beta > 0$. Yet this follows immediately as 
\begin{align*}
|a_0 + \int_0^t \varphi(r) \dot{g}(r) dr| &= \Big|\int_t^{\infty} \varphi(r) \dot{g}(r) dr\Big| \\
& \leq \left(\sup_{s \geq t} |\varphi(s)| \right) \int_t^{\infty} \dot{g}(r) dr \\
& \leq K e^{-\gamma t} [1 - g(t)] \\
& \leq K e^{- (\gamma + \sqrt{2}) t}
\end{align*}
Thus, $\dot{v}$ is exponentially decaying, so that $v(t)$ is bounded and hence
\[
\lim_{t \to \infty} F(t) = \lim_{t \to \infty} v(t) \dot{g}(t) = 0
\]
Moreover, since $v(t)$ is bounded and $\dot{g}(t)$ is exponentially decaying, $F(t)$ is also exponentially decaying by differentiating the equation for $F$.  \qed
\bibliography{Allen_Cahn_Energy_Minimization}{}
\bibliographystyle{plain}
\end{document}